\documentclass[a4paper, 11pt]{article}

\usepackage{graphicx}
\usepackage{lipsum}
\usepackage{amssymb}
\usepackage{amsmath}
\usepackage{amsthm}
\usepackage{enumitem}

\usepackage[utf8]{inputenc}
\usepackage{geometry}
\usepackage{fancyhdr}

\usepackage{titlesec}
\usepackage[colorlinks]{hyperref}

\usepackage{setspace}
\usepackage{dsfont}
\usepackage{xcolor}
\usepackage{xfrac}
\usepackage{tikz-cd}
\usepackage{comment}

\usepackage[blocks]{authblk}
\title{Multi-Split Continuity}

\titleformat{\section}
{\normalfont\LARGE\centering}{\thesection}{1em}{}
\titlespacing{\section}{0cm}{1cm}{0.5cm}

\newtheorem{theorem}{Theorem}[section]
\newtheorem{proposition}[theorem]{Proposition}
\newtheorem{corollary}[theorem]{Corollary}
\newtheorem{conjecture}[theorem]{Conjecture}
\newtheorem{lemma}[theorem]{Lemma}

\theoremstyle{definition} 
\newtheorem*{acknowledgements}{Acknowledgements}
\newtheorem*{contributions}{Contributions}
\newtheorem{definition}[theorem]{Definition}
\newtheorem{remark}[theorem]{Remark}
\newtheorem{example}[theorem]{Example}

\hypersetup{
    colorlinks=true,  
    linkcolor=black,  
    citecolor=black,  
    filecolor=black,
    urlcolor=black,
    linktoc=all       
}

\usepackage{cite}

\begin{document}

\renewcommand*{\Affilfont}{\footnotesize\normalfont}
\renewcommand*{\Authfont}{\Large}

\title{\LARGE \bfseries On Multi-Split Continuity and Split Homeomorphisms}
\author[1]{Finn Michler}
\author[2]{Argha Ghosh\thanks{Corresponding author: argha.ghosh@manipal.edu}}
\affil[1]{D-MATH, ETH Zurich, Rämistrasse 101, CH-8092 Zürich, Switzerland}
\affil[2]{Department of Mathematics, Manipal Institute of Technology, Manipal Academy of Higher Education, Manipal-576104, India}

\date{} 

\maketitle

\vspace{-1cm}

\begin{abstract}
\noindent
    We introduce multi-split continuous functions between topological spaces, a weaker form of continuity that generalizes split continuity while being stable under compositions. We will define the associated star multifunction and pre-multi-split multifunctions. Moreover, we will prove that multi-split continuity naturally emerges as the continuity property of selections of finite usco maps, relating their study to set-valued analysis. Finally, we introduce split homeomorphisms and split homeomorphic spaces, showing that for compact, regular Hausdorff spaces, split homeomorphisms characterize deformations with cuts and subsequent re-glues.
\end{abstract}

\section{Introduction}
\label{Section Introduction}
\noindent
Given a multifunction (set-valued map) $F : X \rightrightarrows Y$, a selection of $F$ is a single-valued map $f : X \to Y$ with $f ( x ) \in F ( x )$ for all $x \in X$. The continuity properties of such selections have been studied extensively. The most prominent examples include the Michael selection theorem on the existence of continuous selections of lower semicontinuous multifunctions \cite{Mi56} and the Castaing representation theorem on the existence of measurable selections of measurable multifunctions \cite{Ca67}. Quasicontinuous selections of usco multifunctions have been the subject of more recent studies in \cite{Ca04, Ho09, Ja22}. In this paper, we will continue this exploration by considering selections of finite upper semicontinuous multifunctions. We will introduce multi-split continuity, which is the natural continuity property of any such selection.

Multi-split continuity is a truly stronger condition compared to Fuller's concept of subcontinuity of functions as introduced in \cite{Fu68} and further studied in \cite{Ho16, Ho12, No05}. A function ${f : X \to Y}$ is called subcontinuous if for every convergent net $\langle x_\lambda \rangle_{\lambda \in \Lambda} \longrightarrow p$ in $X$ the net $\langle f ( x_\lambda ) \rangle_{\lambda \in \Lambda}$ has a cluster point $y_0 \in Y$. Multi-split continuity additionally requires that this cluster point $y_0$ lies in some pre-determined and finite subset $Z_p \subseteq Y$, which we will call a set of extended values of $f$ at $p$. Every multi-split continuous function $f : X \to Y$, where $Y$ is a Hausdorff space, is naturally associated with the star multifunction $f^* : X \rightrightarrows Y, \, p \mapsto Z_p$ which maps any point $p \in X$ to the set of extended values of $f$ at $p$. We will explore the properties of $f^*$ in section \ref{Section The Multisplit Extension}. If the space $Y$ is regular Hausdorff, the star multifunction $f^* : X \rightrightarrows Y$ is upper semicontinuous (Theorem \ref{theorem multi-split extension is usc}), and its graph ${gr ( f^* ) \subseteq X \times Y}$ is the closure of the graph of $f$ (Proposition \ref{proposition graph of f^* is closure of graph of f}).

The concept of multi-split continuity further generalizes the definition of split continuous functions stated by Beer--Bloomfield--Robinson in \cite{Be19} on metric spaces and by the second author in \cite{Gh24} on topological spaces. Split continuous functions have also been the subject of more recent studies in \cite{Ag24, Agg22.2, Agg22.1}. Many of the properties of split continuous functions nicely transfer to multi-split continuous functions. Among others, we will state several equivalent definitions of multi-split continuity using open neighborhoods (Definition \ref{Definition topological multi-split function}, Proposition \ref{proposition equivalent definition only with condition b}), nets (Proposition \ref{Proposition Multi-Split continuity with nets}), and the closures of graphs (Proposition \ref{Proposition equivalent definition using closure of the graph}).

Pre-multi-split multifunctions will be another useful tool in the study of multi-split continuous functions. They are multifunctions, whose selections are multi-split continuous, and we will explore their properties in section \ref{Section Pre-multi-split Functions}. Theorem \ref{theorem finite usc functions are pre-multi-split} states that all non-empty, finite, and upper semicontinuous multifunctions are pre-multi-split, and corollary \ref{corollary union of fin. many multi-split functions is pre-multi-split} enables us to construct new multi-split continuous functions as the selections of finite unions of multi-split continuous functions.

Multi-split continuous functions have one main advantage over split continuous functions, in that they are stable by compositions (Theorem \ref{Theorem compositions of multi-split functions}). This permits us to define split homeomorphisms, multi-split continuous functions with multi-split continuous inverses, and split homeomorphic spaces as an equivalence relation on the class of topological spaces in section \ref{Section Split Homeomorphisms}. This equivalence is truly coarser than homeomorphisms of topological spaces while still preserving some topological properties, such as the compactness of a space (Proposition \ref{Proposition split homeomorphisms and compact spaces}).

Cutting and re-gluing of manifolds has been extensively studied in the field of surgery theory; standard references include \cite{Br12, Ra02}. Surgery, however, only allows for the removal of a tubular neighborhood of an embedded sphere with subsequent re-gluing along its boundary, leading to a self-homeomorphism of the manifold. In section \ref{Section Cuts and Re-glues} we will explore a less restrictive class of deformations, which are well-defined on all compact regular Hausdorff spaces. Given two such spaces $X$ and $Y$, we say that $X$ can be deformed into $Y$ with cuts and subsequent re-glues if there exists a compact space $Z$ and finite-to-one quotient maps $p_X : Z \to X$ and $p_Y : Z \to Y$, together with a right-inverse $p_X^{-1} : X \to Z$ of $p_X$ such that $p_Y \circ p_X^{-1}$ is a bijection. Here, $p_X^{-1}$ and $p_Y$ respectively represent the cutting and the re-gluing maps. We will show that two compact regular Hausdorff spaces are split-homeomorphic if and only if they are equivalent with cuts and subsequent re-glues (Propositions \ref{proposition equivalence with cuts and re-glues implies split homeomorphic} and \ref{proposition split homeomorphic implies equivalence with cuts and subsequent re-glues}). The theory of multi-split continuous functions can therefore be used to study the geometrically interesting procedure of deforming a space together with cuts and subsequent re-glues.

\section{Preliminaries}
\label{Section Preliminaries}
\noindent
Throughout this paper, let $(X , \tau_X )$, $( Y , \tau_Y )$ and $( Z , \tau_Z )$ denote arbitrary topological spaces. We will drop the topologies $\tau_X$, $\tau_Y$, and $\tau_Z$ in notation whenever convenient and impose additional constraints on the spaces $X$, $Y$, and $Z$ wherever necessary. By default, all products of spaces are equipped with the product topology, and all subspaces carry the subspace topology. A subset $V \subseteq X$ is called a \textit{neighborhood of $p$} if there exists an open set $U \subseteq X$ with $p \in U \subseteq V$. We set $\mathfrak{N} ( p )$ to be the set of all neighborhoods of $p$ in $X$. Further, let $cl ( \cdot )$, $int ( \cdot )$, $\partial ( \cdot )$, and $( \cdot )^c$ be the usual \textit{closure, interior, boundary, and complement} operations on subsets of a topological space. Given a metric space $( X' , d )$, a point $p \in X'$ and a constant $\varepsilon > 0$, we set $B_d ( p , \varepsilon)$ to be the open ball with radius $\varepsilon$ and centered at $p$ in $X'$. \\

Let $( X' , d )$ and $( Y' , \rho )$ be two metric spaces, and let $p \in X'$ be a point. A function $f : X' \to Y'$ is called \textit{split continuous at $p$} as introduced by Beer--Bloomfield--Robinson in \cite[Def.~2.1]{Be19} if there exists a $y_p \in Y'$ such that the following two conditions hold.
    \begin{enumerate}[label=\arabic*), topsep=-\parskip, noitemsep]
        \item
        $\forall \varepsilon > 0 \; : \; f \big( B_d ( p , \varepsilon ) \big) \cap B_\rho ( y_p , \varepsilon ) \neq \varnothing $
        \item
        $\forall \varepsilon > 0 \; \exists \delta > 0 \; : \; f \big( B_d ( p , \delta ) \big) \subseteq B_\rho ( f ( p ) , \varepsilon ) \cup B_\rho ( y_p , \varepsilon )$
    \end{enumerate}
Given a function $f : X \to Y$ between topological spaces, the second author introduced the following generalized definition of split continuity in \cite[Def.~3.1]{Gh24}. The function $f$ is called \textit{split continuous at $p \in X$} if there exists a $y_p \in Y$ such that the following two conditions hold.
    \begin{enumerate}[label=\arabic*), topsep=-\parskip, noitemsep]
        \item
        $\forall U \in \mathfrak{N} ( p ) \; \forall V \in \mathfrak{N} ( y_p ) \; : \; f ( U ) \cap V \neq \varnothing $
        \item
        $\forall V_{f(p)} \in \mathfrak{N} \big( f(p) \big) , \, V_{y_p} \in \mathfrak{N} ( y_p ) \; \exists U \in \mathfrak{N} ( p ) \; : \; f ( U ) \subseteq V_{f ( p ) } \cup V_{y_p}$
    \end{enumerate}
The function $f$ is called \textit{globally split continuous} if $f$ is split continuous at every $p \in X$. The two definitions above are equivalent if both spaces $X$ and $Y$ are metrizable. If the space $Y$ is Hausdorff, we can define the \textit{associated star function} $f^* : X \to Y , \, x \mapsto y_x$. We now recite some interesting results about split continuous functions, which we will generalize in this paper. 

Let $f : X \to Y$ be a split continuous function between topological spaces, where $Y$ is a Hausdorff space. Then ${cl \big( gr ( f ) \big) = gr ( f ) \cup gr ( f^* )}$ \cite[Prop.~3.6]{Be19}, \cite[Prop.~3.12]{Gh24}. A function $f : X' \to Y'$ between metric spaces is globally split continuous if and only if there exists a second function $g : X' \to Y'$ such that the multifunction $ \Psi : X' \rightrightarrows Y' , \, x \mapsto \{ f ( x ) , g ( x ) \}$ is upper semicontinuous \cite[Thm.~4.3]{Be19}. A function $f : X \to Y$ with a Hausdorff codomain is continuous if and only if it is split continuous and has a closed graph \cite[Prop.~3.13]{Gh24}. \\

Let $\Lambda$ be an index set together with an upward directed preorder $\preccurlyeq$, that is, for any $\lambda_1 , \lambda_2 \in \Lambda$ there exists some $\lambda_0 \in \Lambda$ with $\lambda_1 , \lambda_2 \preccurlyeq \lambda_0$. A \textit{net} $\langle x_\lambda \rangle_{\lambda \in \Lambda}$ in $X$ is a function $F : \Lambda \to X , \, \lambda \mapsto x_\lambda$. A point $p \in X$ is called a \textit{cluster point of the net} $\langle x_\lambda \rangle_{\lambda \in \Lambda}$ if for every neighborhood $U \subseteq X$ of $p$ and every $\lambda_0 \in \Lambda$ there exists some $\lambda_1 \in \Lambda$ with $\lambda_0 \preccurlyeq \lambda_1$ and $x_{\lambda_1} \in U$. The point $p \in X$ is called a \textit{limit point of the net} $\langle x_\lambda \rangle_{\lambda \in \Lambda}$, and we write $\langle x_\lambda \rangle_{\lambda \in \Lambda} \longrightarrow p$ if for every neighborhood $U \subseteq X$ of $p$ there exists some $\lambda_0 \in \Lambda$ such that for all $\lambda \in \Lambda$ with $\lambda_0 \preccurlyeq \lambda$ we have $x_\lambda \in U$. The net $\langle y_\theta \rangle_{\theta \in \Theta}$ is called a \textit{subnet of the net} $\langle x_\lambda \rangle_{\lambda \in \Lambda}$ if there exists an order-preserving map $\phi : \Theta \to \Lambda$ such that for all $\theta \in \Theta$ we have $x_{\phi ( \theta) } = y_{\theta}$ and $\phi ( \Theta ) \subseteq \Lambda$ is cofinal, that is, for every $\lambda_0 \in \Lambda$ there exists some $\theta_0 \in \Theta$ with $\lambda_0 \preccurlyeq \phi ( \theta_0 )$. We will denote a generic subnet of $\langle x_\lambda \rangle_{\lambda \in \Lambda}$ by $\langle x_{\phi ( \theta )} \rangle_{\theta \in \Theta}$. We now state the following well-known propositions.
\begin{proposition}
\label{proposition space is Hausdorff iff all convergent nets have unique limits}
    A topological space $X$ is Hausdorff if and only if every net $\langle x_\lambda \rangle_{\lambda \in \Lambda}$ in $X$ has at most one limit point \cite[Prop.~1.6.7]{En89}.
\end{proposition}

\begin{proposition}
\label{proposition cluster points of a net are exactly limits of convergent subnets}
    The point $p \in X$ is a cluster point of the net $\langle x_\lambda \rangle_{\lambda \in \Lambda}$ if and only if there exists a convergent subnet $\langle x_{\phi ( \theta ) } \rangle_{\theta \in \Theta} \longrightarrow p$ \cite[Thm.~11.5]{Wi04}. Also, any cluster point of a subnet $\langle x_{\phi ( \theta ) } \rangle_{\theta \in \Theta}$ is also a cluster point of the net $\langle x_\lambda \rangle_{\lambda \in \Lambda}$ \cite[Cor.~11.6]{Wi04}.
\end{proposition}

\begin{proposition}
\label{proposition nets and closure}
    Let $A \subseteq X$ be a subspace. Then $p \in cl ( A )$ if and only if there exists a convergent net $\langle x_\lambda \rangle_{\lambda \in \Lambda} \longrightarrow p$ in $X$, where $x_\lambda \in A$ for all $\lambda \in \Lambda$ \cite[Prop.~1.6.3]{En89}. Further, the set $A$ is closed if and only if every limit point in $X$ of a net in $A$ already lies in $A$ \cite[Cor.~1.6.4]{En89}.
\end{proposition}

\begin{proposition}
\label{proposition nets and compact spaces}
    A space $X$ is compact if and only if every net $\langle x_\lambda \rangle_{\lambda \in \Lambda}$ in $X$ has a cluster point \cite[Thm.~17.4]{Wi04}.
\end{proposition}

A \textit{multifunction} between topological spaces $X$ and $Y$ is a set-valued map \linebreak$F : X \rightrightarrows Y$. We call $F$ \textit{non-empty, open, closed, compact, discrete, or connected} if the image $F ( x )$ possesses the corresponding property for all $x \in X$. For any subset $A \subseteq X$, we set the \textit{image of $A$ under $F$} to be the following set.
$$F ( A ) := \bigcup_{x \in A} F ( x ) \subseteq Y$$
We define the \textit{inverse image $F^{-1} ( B )$ of $B \subseteq Y$ under $F$} and the \textit{core $F^{+1} ( B )$ of $B \subseteq Y$ under $F$} to be the following subsets of $X$.
$$F^{-1} ( B ) := \{ x \in X : F ( x ) \cap B \neq \varnothing \}
\qquad
F^{+1} ( B ) := \{ x \in X : F ( x ) \subseteq B \}$$
We further set the \textit{graph of $F$} to be the following subset of $X \times Y$.
$$gr ( F ) := \{ ( x , y ) \in X \times Y : y \in F ( x ) \} \subseteq X \times Y$$
Given a second multifunction $G : X \rightrightarrows Y$, we define the \textit{union and intersection of $F$ and $G$} to be $F \cup \cap \, G : X \rightrightarrows Y, \, x \mapsto F ( x ) \cup \cap \, G ( x )$. We call $G$ a \textit{submultifunction of $F$} if $gr ( G ) \subseteq gr ( F )$ and call $G$ a \textit{proper submultifunction of $F$} if $gr ( G ) \varsubsetneq gr ( F )$. Let $H : Y \rightrightarrows Z$ be another multifunction. We define the \textit{composition of $H$ and $F$} to be the following multifunction.
$$H \circ F : X \rightrightarrows Z , \quad x \mapsto H \big( F ( x ) \big)$$
Given a single-valued function $f : Y \to X$, we define the \textit{inverse image multifunction} $f^{-1}$ determined by $f$ to be $f^{-1} : X \rightrightarrows Y, \, x \mapsto f^{-1} ( \{ x \} )$. A \textit{selection} of a non-empty multifunction $F : X \rightrightarrows Y$ is a single-valued function $f : X \to Y$ with $f ( x ) \in F ( x )$ for all $x \in X$.

A non-empty multifunction $F : X \rightrightarrows Y$ is called \textit{subcontinuous at $p \in X$} if every selection $f$ of $F$ is subcontinuous at $p \in X$, that is, for every convergent net $\langle x_\lambda \rangle_{\lambda \in \Lambda} \longrightarrow p$, the net $\langle f ( x_\lambda ) \rangle_{\lambda \in \Lambda}$ has a cluster point in $Y$. An arbitrary multifunction $F : X \rightrightarrows Y$ is called \textit{upper semicontinuous at $p \in X$} (or \textit{u.s.c. at $p \in X$} for short) if for every neighborhood $V \subseteq Y$ of $F ( p )$ the core $F^{+1} ( V )$ of $V$ under $F$ is a neighborhood of $p$. We call $F$ \textit{globally upper semicontinuous} if $F$ is u.s.c. at every point $p \in X$. An upper semicontinuous multifunction $F : X \rightrightarrows Y$ with non-empty and compact values is called an \textit{usco map}. $F$ is called \textit{minimal usco} if no proper submultifunction $G$ of $F$ is usco itself. We list the following well-known propositions.

\begin{proposition}
\label{proposition subcontinuous and closed graph iff continuous}
    Let $f : X \to Y$ be a single valued function. If $f$ is subcontinuous and has a closed graph, then $f$ is continuous. If $Y$ is a Hausdorff space, then every continuous function $f : X \to Y$ is subcontinuous and has a closed graph \cite[Thm.~3.4]{Fu68}.
\end{proposition}

\begin{proposition}
\label{proposition connection of upper semicontinuity and subcontinuity}
    Let $F : X \rightrightarrows Y$ be a multifunction and let $p \in X$ be a point. If $F$ is upper semicontinuous at $p$ and $F ( p ) \subseteq Y$ is compact, then $F$ is subcontinuous at $p$ \cite[Thm.~2.5]{Le78}. If $F$ is subcontinuous and has a closed graph, then $F$ is upper semicontinuous \cite[Thm.~3.1]{Sm75}.
\end{proposition}

\begin{proposition}
\label{Proposition usco maps send compact sets to compact sets}
    Let $F : X \rightrightarrows Y$ be an usco map. If $K \subseteq X$ is compact, then so is $F ( K ) \subseteq Y$ \cite[Prop.~1.1.7]{Ho21}.
\end{proposition}

\begin{proposition}
\label{proposition closed graph theorem for usc multifunctions}
    Let $F : X \rightrightarrows Y$ be an usco map. If $Y$ is Hausdorff, then the graph $gr ( F )$ of $F$ is closed in $X \times Y$ \cite[Prop.~1.1.8]{Ho21}. Now let $Y$ be compact and $G : X \rightrightarrows Y$ be a multifunction such that the set $\{ x \in X : G ( x ) \neq \varnothing \} \subseteq X$ is dense in $X$. If $G$ has a closed graph, then $G$ is an usco map \cite[Cor.~1.1.15]{Ho21}.
\end{proposition}

\begin{proposition}
\label{proposition inverse image multifunction in u.s.c. iff f is closed}
    Let $f : X \to Y$ be a continuous map, and let $f^{-1} : X \rightrightarrows Y$ be the inverse image multifunction determined by $f$. Then $f^{-1}$ is globally upper semicontinuous if and only if $f$ is closed \cite[Sec.~6.2, Ex.~4a]{Be93}.
\end{proposition}

\section{Multi-Split Continuity}
\label{Section Multi-Split Continuity}
\noindent
In this section, we will introduce the notion of multi-split continuous functions between topological spaces. We will define sets of extended values and prove some basic properties of multi-split continuous functions. Theorem \ref{Theorem compositions of multi-split functions}, describing the stability of multi-split continuous functions by compositions, as well as proposition \ref{Proposition Multi-Split continuity with nets}, the characterization of multi-split continuity using nets, will be of great value throughout this paper.

\begin{definition}
\label{Definition topological multi-split function}
    Let $f : X \to Y$ be a function and $p \in X$ be a point. We call $f$ \textit{multi-split continuous at $p$} if there exists a finite non-empty collection of distinct values $Z_p := \{ y_1 , ... \, , y_n \} \subseteq Y$ such that the following two conditions hold:
    \begin{enumerate}[label=\alph*), topsep=-\parskip, noitemsep]
        \item
        \label{definition topological multi-split continuity a}
        For every $y_i \in Z_p$ and for all neighborhoods $U \subseteq X$ of $p$ and $V_i \subseteq Y$ of $y_i$ we have $f ( U ) \cap V_i \neq \varnothing$.
        \item
        \label{definition topological multi-split continuity b}
        For every neighborhood $V \subseteq Y$ of $Z_p$ there exists a neighborhood $U \subseteq X$ of $p$ such that $f ( U ) \subseteq V$.
    \end{enumerate}
    We call $Z_p$ a \textit{set of extended values of $f$ at $p$}. If $f$ is multi-split continuous at $p$ for every $p \in X$, we say that $f$ is \textit{(globally) multi-split continuous}. We set $mSplit ( X , Y )$ to be the set of all globally multi-split continuous functions $f : X \to Y$.
\end{definition}

The concept of multi-split continuity generalizes split continuity of functions as introduced by Beer--Bloomfield--Robinson on metric spaces in \cite{Be19} and by the second author on topological spaces in \cite{Gh24}. Some of the following most basic properties and their proofs are inspired by their work.

\begin{remark}
\label{Remark multi-split versus continuous and split-continuous}
    Multi-split continuity of the function $f : X \to Y$ at a point \linebreak$p \in X$ naturally reduces to the usual continuity of functions between topological spaces in the case $Z_p = \{ f ( p ) \}$ and to split continuity of functions as introduced in \cite{Gh24} in the case $| Z_p | = 2$ with $f ( p ) \in Z_p$.
\end{remark}

Condition \ref{Definition topological multi-split function}\ref{definition topological multi-split continuity a} is redundant to a certain extent, as the following proposition shows. We will, however, stick with the definition above. The set of extended values $Z_p$ of a given multi-split continuous function $f$ at $p$ will be of great interest to us later, and including it in the definition will be very convenient.

\begin{proposition}
\label{proposition equivalent definition only with condition b}
    Let $f : X \to Y$ be a function and let $p \in X$ be a point. Then the following are equivalent.
    \begin{enumerate}[label=\arabic*), nosep]
        \item
        \label{proposition equivalent definition only b 1}
        $f$ is multi-split continuous at $p$.
        \item 
        \label{proposition equivalent definition only b 2}
        There exists a non-empty finite subset $Z \subseteq Y$ such that condition \ref{Definition topological multi-split function}\ref{definition topological multi-split continuity b} holds for the set $Z$. That is, for every neighborhood $V \subseteq Y$ of $Z$ there exists a neighborhood $U \subseteq X$ of $p$ such that $f ( U ) \subseteq V$.
    \end{enumerate}
\end{proposition}

\begin{proof}
    Statement \ref{proposition equivalent definition only b 1} clearly implies statement \ref{proposition equivalent definition only b 2} by choosing $Z := Z_p$. We verify the second implication. Let $Z \subseteq Y$ be a non-empty finite subset such that statement \ref{proposition equivalent definition only b 2} holds. We define $Z_p \subseteq Z$ as follows.
    \begin{equation*}
        Z_p := \{ z \in Z : \text{ condition \ref{Definition topological multi-split function}\ref{definition topological multi-split continuity a} holds for } f \text{ at } z \}
    \end{equation*}
    We show that $Z_p$ is a set of extended values of $f$ at $p$. $Z_p$ is a subset of $Z$, which is finite. Therefore, $Z_p$ is a finite set itself. Assume, for contradiction, that $Z_p$ is empty. In particular, for each $y_i \in Z$ there exists a neighborhood $V_i \subseteq Y$ of $y_i$ such that $f ( p ) \notin V_i$. The set $V := \bigcup_{y_i \in Z} V_i \subseteq Y$ then is a neighborhood of $Z$ with $f ( p ) \notin V$. This is a contradiction, since $f ( p ) \in f ( U )$ for every neighborhood $U \subseteq X$ of $p$. Thus, $Z_p$ is non-empty.
    
    It remains to verify conditions \ref{Definition topological multi-split function}\ref{definition topological multi-split continuity a} and \ref{Definition topological multi-split function}\ref{definition topological multi-split continuity b} for the set $Z_p$. Condition \ref{Definition topological multi-split function}\ref{definition topological multi-split continuity a} holds trivially. Thus, let $V \subseteq Y$ be a neighborhood of $Z_p$. Let \linebreak$Z \backslash Z_p =: \{ y_1 , ... \, , y_r \}$. By negating condition \ref{Definition topological multi-split function}\ref{definition topological multi-split continuity a} for each $y_i$ with ${i = 1 , ... \, , r}$, we can find neighborhoods $U_i \subseteq X$ of $p$ and $V_i \subseteq Y$ of $y_i$ such that \linebreak${f ( U_i ) \cap V_i = \varnothing}$. We define the set $V'$ as follows.
    $$V' \, := \, V \cup \bigcup_{i=1}^r V_i \, \subseteq \, Y$$
    $V'$ then is a neighborhood of $Z$ in $Y$. By statement \ref{proposition equivalent definition only b 2}, there exists a neighborhood $U' \subseteq X$ of $p$ with $f ( U' ) \subseteq V'$. We set $U \subseteq X$ to be the following set.
    $$U \, := \, U' \cap \bigcap_{i = 1}^r U_i \, \subseteq \, X$$
    Note that $U$ is a neighborhood of $p$ in $X$ with $f ( U ) \cap V_i \subseteq f ( U_i ) \cap V_i = \varnothing$ for each $i = 1 , ... \, , r$. Thus, $f ( U ) \subseteq V$ and $f$ is multi-split continuous at $p$ as required. This concludes the second implication and therefore the proof.
\end{proof}

This proof implies that the set $Z$ specified in condition \ref{proposition equivalent definition only b 2} in proposition \ref{proposition equivalent definition only with condition b} above contains a set of extended values $Z_p$ of $f$ at $p$. This set of extended values might, however, not be unique. Also, $f ( p )$ might not be an element of $Z_p$. Using the following proposition, we may still choose a set of extended values $Z_p$ of $f$ at $p$ with $f ( p ) \in Z_p$.

\begin{proposition}
\label{proposition Z_p contains f(p)}
    Given a function $f : X \to Y$ and a point $p \in X$ such that $f$ is multi-split continuous at $p$, we can find a set of extended values $Z_p$ of $f$ at $p$ with $f ( x ) \in Z_p$.
\end{proposition}

\begin{proof}
    Let $Z_p$ be a set of extended values of $f$ at $p$. We will show that $\widetilde{Z}_p := Z_p \cup \{ f ( p ) \}$ also fulfills conditions \ref{Definition topological multi-split function}\ref{definition topological multi-split continuity a} and \ref{Definition topological multi-split function}\ref{definition topological multi-split continuity b}. Condition \ref{Definition topological multi-split function}\ref{definition topological multi-split continuity a} holds trivially for $f ( p ) \in \widetilde{Z}_p$. Further note that any open neighborhood $V \subseteq Y$ of $\widetilde{Z}_p$ also is an open neighborhood of $Z_p \subseteq \widetilde{Z}_p$ in $Y$. Thus condition \ref{Definition topological multi-split function}\ref{definition topological multi-split continuity b} for the set $Z_p$ directly implies condition \ref{Definition topological multi-split function}\ref{definition topological multi-split continuity b} for the set $\widetilde{Z}_p$.
\end{proof}

The following lemma further describes the possible sets of extended values of a function $f$ at a given point $p$.

\begin{lemma}
\label{lemma unions of sets of extended values}
    Let $f : X \to Y$ be multi-split continuous at $p \in X$ and let $Z_p , \widetilde{Z}_p \subseteq Y$ be sets of extended values of $f$ at $p$. Then every set $B \subseteq Y$ with $Z_p \subseteq B \subseteq Z_p \cup \widetilde{Z}_p$ also is a set of extended values of $f$ at $p$.
\end{lemma}

\begin{proof}
    Let $B \subseteq Y$ be a set as in the statement above. $B$ then contains the non-empty set $Z_p$ and is therefore itself non-empty. Condition \ref{Definition topological multi-split function}\ref{definition topological multi-split continuity a} holds for all $y \in Z_p \cup \widetilde{Z}_p$ and therefore also for all $y \in B$. Lastly note that every neighborhood $V \subseteq Y$ of $B$ also is a neighborhood of $Z_p$. Condition \ref{Definition topological multi-split function}\ref{definition topological multi-split continuity b} for the set $Z_p$ therefore directly implies condition \ref{Definition topological multi-split function}\ref{definition topological multi-split continuity b} for the set $B$. Thus, $B$ also is a set of extended values of $f$ at $p$.
\end{proof}

If we want the set of extended values $Z_p$ of $f$ at $p$ to be unique, we have to assume that the space $Y$ is Hausdorff. Recall the following fact about a Hausdorff space $Y$. For every finite set $W := \{w_1 , ... \, , w_n \} \subseteq Y$ of $n$ distinct points, there exists a collection of pairwise disjoint neighborhoods $\{V_1 , ... \, , V_n \}$, where for each $i = 1 , ... \, , n$, the set $V_i \subseteq Y$ is a neighborhood of $w_i$.

The following proof generalizes an argument by the second author in \cite{Gh24}.

\begin{proposition}
\label{Proposition uniqueness of set of extended values}
    Let $Y$ be a Hausdorff space and $f : X \to Y$ a function that is multi-split continuous at $p \in X$. Then the set of extended values $Z_p$ of $f$ at $p$ is unique.
\end{proposition}

\begin{proof}
    Assume, for contradiction, that both $Z_p := \{ x_1 , ... \, , x_n \}$ and \linebreak$\widetilde{Z}_p := \{x_0, y_1 , ... \, , y_m \}$ are sets of extended values of $f$ at $p$ with $Z_p \neq \widetilde{Z}_p$. WLOG, let $x_0 \in \widetilde{Z}_p \backslash Z_p$. By lemma \ref{lemma unions of sets of extended values} the set $Z_p \sqcup \{ x_0 \}$ is also a set of extended values of $f$ at $p$. We may therefore WLOG assume that $\widetilde{Z}_p = Z_p \sqcup \{ x_0 \}$.

    Using that $Y$ is a Hausdorff space, for $i = 0 , ... \, , n$ we choose pairwise disjoint neighborhoods $V_i \subseteq Y$ of $x_i$. The set $V := \bigcup_{i=1}^n V_i$ then is a neighborhood of $Z_p$ in $Y$, and by condition \ref{Definition topological multi-split function}\ref{definition topological multi-split continuity b} of the multi-split continuity of $f$ at $p$, there exists a neighborhood $U\subseteq X$ of $p$ with $f ( U ) \subseteq V$. In particular, $f ( U ) \cap V_{0} = \varnothing$, since $V_0 \cap \big( \bigcup_{i=1}^n V_i \big) = \varnothing$. This contradicts condition \ref{Definition topological multi-split function}\ref{definition topological multi-split continuity a} at the point $x_0 \in \widetilde{Z}_p$. Therefore the set of extended values $Z_p$ of $f$ at $p$ is uniquely determined.
\end{proof}

The following theorem is one main motivation for introducing this new but somewhat tedious notion of continuity.

\begin{theorem}
\label{Theorem compositions of multi-split functions}
    Given two globally multi-split continuous functions $f, g$ with $f \in mSplit( X , Y )$ and ${g \in mSplit ( Y , Z ) }$, their composition is also globally multi-split continuous. That is $g \circ f \in mSplit ( X , Z )$.
\end{theorem}

\begin{proof}
    Let $p \in X$ be arbitrary and let $Z_p^f \subseteq Y$ be a set of extended values of $f$ at $p$. Further, for every $y \in Z_p^f$, let $Z_y^g \subseteq Z$ be a set of extended values of $g$ at $y$. We define $\widetilde{Z} \subseteq Z$ as follows.
    $$\widetilde{Z} := \bigcup_{y \in Z_p^f} Z_y^g$$
    We verify that the set $\widetilde{Z}$ fulfills condition \ref{proposition equivalent definition only b 2} from proposition \ref{proposition equivalent definition only with condition b}. $\widetilde{Z}$ is the union of finite sets and thus itself finite. Let $V \subseteq Z$ be an arbitrary neighborhood of $\widetilde{Z}$. Let $y \in Z_p^f$ be arbitrary. $V$ then also is a neighborhood in $Z$ of the set of extended values $Z_y^g$ of $g$ at $y$. By the multi-split continuity of $g$ at $y$, there exists a neighborhood $W_y \subseteq Y$ of $y$ such that $g ( W_y ) \subseteq V$. We define the set $W \subseteq Y$ as follows.
    $$W := \bigcup_{y \in Z_p^f} W_y$$
    $W$ then is a neighborhood in $Y$ of the set of extended values $Z_p^f$ of $f$ at $p$. By the multi-split continuity of $f$ at $p$, there exists a neighborhood $U \subseteq X$ of $p$ with $f ( U ) \subseteq W$. We therefore have $( g \circ f ) ( U ) \subseteq V$ as required and may conclude the proof using proposition \ref{proposition equivalent definition only with condition b}.
\end{proof}

We can characterize the multi-split continuity of $f : X \to Y$ for arbitrary spaces $X$ and $Y$ using nets. The following proposition is a generalization of a theorem by the second author and this proof adapts the argument presented in \cite[Thm.~3.3]{Gh24}.

\begin{proposition}
\label{Proposition Multi-Split continuity with nets}
    Let $f : X \to Y$ be a function, and let $p \in X$ be a point. Then the following are equivalent:
    \begin{enumerate}[label=\arabic*), left=0pt, nosep]
        \item
        $f$ is multi-split continuous at $p$
        \item
        \label{Definition sequential multi-split continuity}
        There exists a finite non-empty collection $Z_p := \{ y_1 , ... \, , y_n \} \subseteq Y$ such that the following two conditions hold:
        \begin{enumerate}[label=\alph*'), nosep]
            \item
            \label{definition sequential multi-split continuity a}
            For every $y_i \in Z_p$ there exists a convergent net $\langle x_\lambda \rangle_{\lambda \in \Lambda} \longrightarrow p$ in $X$ with $\langle f ( x_\lambda ) \rangle_{\lambda \in \Lambda} \longrightarrow y_i$ in $Y$.
            \item
            \label{definition sequential multi-split continuity b}
            For every convergent net $ \langle x_\lambda \rangle_{\lambda \in \Lambda} \longrightarrow p$ in $X$, there exists some \linebreak$y_i \in Z_p$ which is a cluster point of the net $\langle f ( x_\lambda ) \rangle_{\lambda \in \Lambda}$ in $Y$.
        \end{enumerate}
    \end{enumerate}
\end{proposition}

\begin{proof}
    We show two implications.
    \begin{description}[font=\normalfont\itshape, leftmargin=15pt, nosep, listparindent=\parindent] 
        \item[$1) \Rightarrow 2) :$]
         Let $f : X \to Y$ be multi-split continuous at $p$ and let \linebreak$Z_p := \{ y_1 , ... \, , y_n \}$ be a set of extended values of $f$ at $p$. Let $y_i \in Z_p$ be arbitrary. We set $\Lambda_i$ to be the set of all pairs of neighborhoods $( U , V )$, where $U \subseteq X$ is a neighborhood of $p$ and $V \subseteq Y$ is a neighborhood of $y_i$. We define the preorder $\preccurlyeq$ on $\Lambda_i$ as follows.
         $$ ( U , V ) \preccurlyeq ( \widetilde{U} , \widetilde{V} ) \quad : \Longleftrightarrow \quad \widetilde{U} \subseteq U \text{ and } \widetilde{V} \subseteq V$$
         Note that the intersection of two neighborhoods of $p$ in $X$ (of $y_i$ in $Y$) is itself a neighborhood of $p$ in $X$ (of $y_i$ in $Y$). Given two pairs of neighborhoods $( U , V ) , ( U' , V' ) \in \Lambda_i$ we therefore have $(\widetilde{U} , \widetilde{V} ) := ( U \cap U' , V \cap V' ) \in \Lambda_i$ with $( U , V ) , ( U' , V' ) \preccurlyeq ( \widetilde{U} , \widetilde{V} )$ and $\preccurlyeq$ is an upward directed preorder. We may therefore use $\Lambda_i$ as the index set for a net.
         
         By condition \ref{Definition topological multi-split function}\ref{definition topological multi-split continuity a}, for every pair of neighborhoods $( U , V ) \in \Lambda_i$, we can choose some $x_U^V \in U$ with $f ( x_U^V ) \in V$. For any neighborhood $U$ of $p$ in $X$, and $( U' , V' ) \in \Lambda_i$ with $( U , Y ) \preccurlyeq ( U' , V' )$, we have $x_{U'}^{V'} \in U' \subseteq U$. The net $\langle x_U^V \rangle_{( U , V ) \in \Lambda_i }$ therefore converges to $p$ in $X$. Similarly note that net $\langle f ( x_U^V ) \rangle_{( U , V ) \in \Lambda_i}$ converges to $y_i$ in $Y$. This verifies condition \ref{Proposition Multi-Split continuity with nets}\ref{definition sequential multi-split continuity a}.

        For contradiction, now assume that condition \ref{Proposition Multi-Split continuity with nets}\ref{definition sequential multi-split continuity b} is false. That is, there exists a convergent net $\langle x_\lambda \rangle_{\lambda \in \Lambda} \longrightarrow p$ in $X$ such that the net $\langle f ( x_\lambda ) \rangle_{\lambda \in \Lambda} $ has no cluster point $y_i \in Z_p$. Therefore, there exists a neighborhood $V \subseteq Y$ of $Z_p$ and a lower bound $\lambda_0 \in \Lambda$, such that for all $\lambda \in \Lambda$ with $\lambda_0 \preccurlyeq \lambda$, we have $f ( x_\lambda ) \notin V$. By condition \ref{Definition topological multi-split function}\ref{definition topological multi-split continuity b}, there exists a neighborhood $U \subseteq X$ of $p$, such that $f ( U ) \subseteq V$. By convergence of the net $\langle x_\lambda \rangle_{\lambda \in \Lambda}$, we can find some $\lambda_1 \in \Lambda$ with $\lambda_0 \preccurlyeq \lambda_1$ and $x_{\lambda_1} \in U$. We therefore have $f ( x_{ \lambda_1 } ) \in V$, a contradiction. Thus, condition \ref{Proposition Multi-Split continuity with nets}\ref{definition sequential multi-split continuity b} must hold as well.
        \item[$2) \Rightarrow 1) :$]
        By definition, the set $Z_p$ specified in statement \ref{Definition sequential multi-split continuity} above is non-empty and finite. We show that the set $Z_p$ also fulfills conditions \ref{Definition topological multi-split function}\ref{definition topological multi-split continuity a} and \ref{Definition topological multi-split function}\ref{definition topological multi-split continuity b} from the definition of a multi-split continuous function. Let $U\subseteq X$ be a neighborhood of $p$ and $V \subseteq Y$ be a neighborhood of $y_i \in Z_p$. By statement \ref{Proposition Multi-Split continuity with nets}\ref{definition sequential multi-split continuity a}, we can find a net $\langle x_\lambda \rangle_{\lambda \in \Lambda}$ in $X$ convergent to $p$ such that the net $\langle f ( x_\lambda ) \rangle_{\lambda \in \Lambda}$ converges to $y_i$ in $Y$. By the definition of convergence of nets in topological spaces, there exist indices $\lambda'_0, \widetilde{\lambda}_0 \in \Lambda$ such that for all $\lambda' \in \Lambda$ with $\lambda'_0 \preccurlyeq \lambda'$ we have $x_{\lambda'} \in U$ and for all $\widetilde{\lambda} \in \Lambda$ with $\widetilde{\lambda}_0 \preccurlyeq \widetilde{\lambda}$ we have $f ( x_{\widetilde{\lambda}} ) \in V$. Taking a common upper bound $\lambda_0 \in \Lambda$, that is, $\lambda'_0 , \widetilde{\lambda}_0 \preccurlyeq \lambda_0$, we get $x_{\lambda_0} \in U$ and $f ( x_{\lambda_0} ) \in V$ and thus $f( U ) \cap V \neq \varnothing$, verifying condition \ref{Definition topological multi-split function}\ref{definition topological multi-split continuity a}.

        Now let $V \subseteq Y$ be a neighborhood of $Z_p$. The set $\Lambda$ of all neighborhoods of $p$ in $X$ is partially ordered and upward directed by inclusion. We may therefore use it as the index set for a net. For contradiction, assume that condition \ref{Definition topological multi-split function}\ref{definition topological multi-split continuity b} does not hold. That is, assume that for every $U \in \Lambda$ we have $f ( U ) \cap V^c \neq \varnothing$. Now choose some $x_U \in U$ with $f ( x_U ) \notin V$ for every $U \in \Lambda$. The net $\langle x_U \rangle_{U \in \Lambda}$ then clearly converges to $p$ in $X$. However, none of the $f ( x_U )$ are elements of $V$, which is a neighborhood of $Z_p$. The net $\langle f ( x_U ) \rangle_{U \in \Lambda}$ therefore has no cluster points in $Z_p$. This contradicts condition \ref{Proposition Multi-Split continuity with nets}\ref{definition sequential multi-split continuity b}, and condition \ref{Definition topological multi-split function}\ref{definition topological multi-split continuity b} must hold.
    \end{description}
This concludes the second implication and therefore the proof.
\end{proof}

Proposition \ref{Proposition Multi-Split continuity with nets} now leads to the following more explicit expression for the set of extended values $Z_p$ of $f$ at $p$ using convergent nets.

\begin{proposition}
\label{proposition net of function values can only have fin. many cluster points}
    Let $Y$ be a Hausdorff space and let $f : X \to Y$ be multi-split continuous at $p \in X$. Then the set $X_p$, defined below explicitly, is the set of extended values of $f$ at $p$.
    \begin{equation*}
        X_p := \left\{
        y \in Y \;\middle|\;
        \parbox{0.65\textwidth}{
            there exists a convergent net $\langle x_\lambda \rangle_{\lambda \in \Lambda} \longrightarrow p$ in $X$ such that $y$ is a cluster point of the net $\langle f ( x_\lambda ) \rangle_{\lambda \in \Lambda}$ }
        \right\}
    \end{equation*}
\end{proposition}

\begin{proof}
    Let $Z_p$ be the set of extended values of $f$ at $p$. We show two inclusions. Condition \ref{Proposition Multi-Split continuity with nets}\ref{definition sequential multi-split continuity a} of the characterization of multi-split continuity using nets clearly implies $Z_p \subseteq X_p$. Now let $\langle x_\lambda \rangle_{\lambda \in \Lambda} \longrightarrow p$ be a convergent net in $X$ and let $y \in Y$ be a cluster point of the net $\langle f ( x_\lambda ) \rangle_{\lambda \in \Lambda}$ in $Y$. Using proposition \ref{proposition cluster points of a net are exactly limits of convergent subnets}, we can find a subnet $\langle x_{ \phi ( \theta ) } \rangle_{ \theta \in \Theta }$ such that $y$ is a limit point of the net $\langle f ( x_{ \phi ( \theta ) } ) \rangle_{\theta \in \Theta}$. $Y$ is Hausdorff, and thus $y \in Y$ is the unique limit point and therefore also the unique cluster point of this net by proposition \ref{proposition space is Hausdorff iff all convergent nets have unique limits}. Condition \ref{Proposition Multi-Split continuity with nets}\ref{definition sequential multi-split continuity b} now requires that this cluster point $y$ has to lie in $Z_p$ and therefore $Z_p \supseteq X_p$.
\end{proof}

Recall the following fact from general topology. Let $A \subseteq X$ be a subset of a topological space and let $p \in X$ be a point. We then have $p \in cl ( A )$ if and only if for every neighborhood $V \subseteq X$ of $p$ we have $V \cap A \neq \varnothing$ \cite[Prop.~1.1.1]{En89}.

\begin{lemma}
\label{Lemma theorem usc of f^*}
    Let $Y$ be a Hausdorff space, let $f : X \to Y$ be multi-split continuous at $p \in X$, and let $Z_p \subseteq Y$ be the set of extended values of $f$ at $p$. For every neighborhood $U \subseteq X$ of $p$ we then have that $Z_p \subseteq cl \big( f ( U ) \big)$.
\end{lemma}

\begin{proof}
    Let $y_i \in Z_p$ be an arbitrary extended value of $f$ at $p$. Let $V \subseteq Y$ be a neighborhood of $y_i$. By condition \ref{Definition topological multi-split function}\ref{definition topological multi-split continuity a} in the definition of multi-split continuity of $f$ at $p$ we have $f ( U ) \cap V \neq \varnothing$. Thus, every neighborhood of $y_i$ intersects $f ( U )$. By the fact quoted above we then have $y_i \in cl \big( f ( U ) \big)$ and therefore $Z_p \subseteq cl \big( f ( U ) \big)$.
\end{proof}

\begin{proposition}
\label{proposition explicit formula for set of extended values}
    Let $Y$ be a Hausdorff space and let $f : X \to Y$ be multi-split continuous at $p \in X$. Then the set $X_p$, defined explicitly as follows, is the set of extended values of $f$ at $p$.
    \begin{equation*}
        X_p := 
        \bigcap_{U \in \, \mathfrak{N} ( p )} cl \big( f ( U ) \big)
    \end{equation*}
\end{proposition}

\begin{proof}
    Let $Z_p \subseteq Y$ be the set of extended values of $f$ at $p$. We show two inclusions. Lemma \ref{Lemma theorem usc of f^*} immediately implies $Z_p \subseteq X_p$. Now let $y_i \in X_p$ be arbitrary. We show that $Z_p \backslash \{ y_i \}$ fails condition \ref{Proposition Multi-Split continuity with nets}\ref{definition sequential multi-split continuity b} of the characterization of multi-split continuity using nets.

    Let $\Lambda_i$ be the set of all pairs of neighborhoods $( U , V )$, where $U \subseteq X$ is a neighborhood of $p$ and $V \subseteq Y$ is a neighborhood of $y_i$. We define the upward directed pre-order $\preccurlyeq$ on $\Lambda_i$ as follows.
    $$ ( U , V ) \preccurlyeq ( \widetilde{U} , \widetilde{V} ) \quad : \Longleftrightarrow \quad \widetilde{U} \subseteq U \text{ and } \widetilde{V} \subseteq V$$
    Let $( U , V ) \in \Lambda_i$ be an arbitrary pair of neighborhoods. For each $( U , V ) \in \Lambda_i$ we can choose some $x_V^U \in U$ with $f ( x_V^U ) \in V$ by the definition of $X_p$. We then clearly have the following two convergent nets.
    $$\langle x_V^U \rangle_{( U , V ) \in \Lambda_i} \longrightarrow p
    \qquad \text{and} \qquad
    \langle f ( x_V^U ) \rangle_{( U , V ) \in \Lambda_i} \longrightarrow y_i$$
    The space $Y$ was assumed to be Hausdorff. Using proposition \ref{proposition space is Hausdorff iff all convergent nets have unique limits}, we may conclude that $y_i$ is the only cluster point of the net $\langle f ( x_V^U ) \rangle_{( U , V ) \in \Lambda_i}$. Condition \ref{Proposition Multi-Split continuity with nets}\ref{definition sequential multi-split continuity b} therefore fails for the set $Z_p \backslash \{ y_i \}$ and we get $X_p \subseteq Z_p$. This concludes the second inclusion and therefore the proof.
\end{proof}

\section{The Star Multifunction}
\label{Section The Multisplit Extension}
\noindent
Any globally multi-split continuous function $f : X \to Y$, where $Y$ is Hausdorff, is naturally associated with its star multifunction $f^* : X \rightrightarrows Y$. In this section we will explore the properties of $f^*$. Theorem \ref{theorem multi-split extension is usc} establishes that $f^*$ is upper semicontinuous, linking multi-split continuous functions to usco maps. Proposition \ref{proposition graph of f^* is closure of graph of f} reveals the intuitive geometric relationship of $f$ and $f^*$, stating that $cl \big( gr ( f ) \big) = gr ( f^* )$.

\begin{definition}
\label{definition multi-split extension}
    Let $Y$ be a Hausdorff space and let $f \in mSplit ( X , Y )$. For every $p \in X$, let $Z_p \subseteq Y$ be the set of extended values of $f$ at $p$. We define the corresponding \textit{star multifunction $f^*$} to be the following.
    $$f^* : X \rightrightarrows Y , \quad p \mapsto Z_p$$
\end{definition}

Note that the set $Z_p \subseteq Y$ is unique for all $p \in X$ by proposition \ref{Proposition uniqueness of set of extended values}. Thus, this definition is well-posed. Further, note that this definition of the star multifunction does not match the definition given by Beer--Bloomfield--Robinson in \cite{Be19}. They defined the associated star function to be single-valued by excluding the function value $f ( p )$ from $Z_p$ at points of strict split continuity. We will turn $f^*$ into a single-valued function in section \ref{Section Pre-multi-split Functions} by considering selections but will stick to this more general definition for now.

Recall the following well-known facts from general topology. Let $X$ be a topological space and let $A_1 , ... \, , A_n \subseteq X$ be an arbitrary finite collection of subsets. We then have $cl ( \bigcup_{i=1}^n A_i ) = \bigcup_{i=1}^n cl ( A_i )$ \cite[Thm.~1.1.3.~CO3]{En89}.
If $Y$ is a regular Hausdorff space, then, for every $y \in Y$ and every neighborhood $V \subseteq Y$ of $y$, there exists a neighborhood $U \subseteq Y$ of $y$ with $cl ( U ) \subseteq V$ \cite[Prop.~1.5.5]{En89}.

\begin{lemma}
\label{lemma theorem usc f^* 2}
    Given a regular Hausdorff space $Y$ and a finite subset \linebreak$Z := \{ y_1 , ... \, , y_n \} \subseteq Y$ as well as a neighborhood $V \subseteq Y$ of $Z$, there exists a neighborhood $A \subseteq Y$ of $Z$ with $cl ( A ) \subseteq V$.
\end{lemma}

\begin{proof}
    Let $V \subseteq Y$ be a neighborhood of $Z$. For every $y_i \in Z$ we can then find a neighborhood $A_i \subseteq Y$ of $y_i$ with $cl ( A_i ) \subseteq V$ for $i = 1 , ... \, , n$ using the regularity of $Y$. Then $A := \bigcup_{i=1}^n A_i$ is a neighborhood of $Z$ in $Y$ with $cl ( A ) = \bigcup_{i=1}^n cl ( A_i ) \subseteq V$.
\end{proof}

The following theorem generalizes a statement by Beer--Bloomfield--{\linebreak}Robinson in \cite[Prop.~4.2]{Be19}, and the proof presented here extends their ideas.

\begin{theorem}
\label{theorem multi-split extension is usc}
    Let $Y$ be a regular Hausdorff space and $f \in mSplit ( X , Y )$. Then the star multifunction $f^* : X \rightrightarrows Y$ is globally upper semicontinuous.
\end{theorem}

\begin{proof}
    Let $p \in X$ be arbitrary and $Z_p := f^* ( p )$ be the set of extended values of $f$ at $p$. Let $V \subseteq Y$ be an arbitrary neighborhood of $Z_p$. Using lemma \ref{lemma theorem usc f^* 2}, we can find a neighborhood $A \subseteq Y$ of $Z_p$ with $cl ( A ) \subseteq V$. By condition \ref{Definition topological multi-split function}\ref{definition topological multi-split continuity b} of the definition of the multi-split continuity of $f$ at $p$, there exists a neighborhood $U \subseteq X$ of $p$ with $f ( U ) \subseteq A$. We now have $Z_p \subseteq cl ( A ) \subseteq V$ by lemma \ref{Lemma theorem usc of f^*} as required and may conclude that $f^*$ is upper semicontinuous at every $p \in X$.
\end{proof}

Finite subsets of topological spaces are compact. The following corollary is an immediate consequence.

\begin{corollary}
\label{corollary f^* is usco map}
    Let $f \in mSplit( X , Y )$, where $Y$ is a regular Hausdorff space. Then the star multifunction $f^* : X \rightrightarrows Y$ is an usco map.
\end{corollary}

As an immediate consequence of this corollary and proposition \ref{Proposition usco maps send compact sets to compact sets}, we can state the following result.

\begin{corollary}
\label{corollary f^* sends cpt sets to cpt sets}
    Let $Y$ be a regular Hausdorff space and let $f \in mSplit ( X , Y )$. Then for every compact set $K \subseteq X$, the set $f^* ( K ) \subseteq Y$ also is compact.
\end{corollary}

The next corollary follows immediately by combining proposition \ref{proposition closed graph theorem for usc multifunctions} and corollary \ref{corollary f^* is usco map}.

\begin{corollary}
\label{corollary f^* has closed graph}
    Let $Y$ be regular Hausdorff and let $f \in mSplit ( X , Y )$. Then $f^* : X \rightrightarrows Y$ has a closed graph.
\end{corollary}

We now generalize another proposition by Beer--Bloomfield--Robinson in \cite[Prop.~3.6]{Be19} to the realm of multi-split continuous functions.

\begin{proposition}
\label{proposition graph of f^* is closure of graph of f}
    Let $Y$ be regular Hausdorff and let $f \in mSplit ( X , Y )$. Then $cl \big( gr ( f ) \big) = gr ( f^* ) \subseteq X \times Y$.
\end{proposition}

\begin{proof}
    We verify two inclusions.
    \begin{description}[font=\normalfont\itshape, leftmargin=15pt, nosep]
        \item[$\subseteq :$]
        For every $x \in X$ we have $f ( x ) \in f^* ( x )$ and thus also that $gr ( f ) \subseteq gr ( f^* )$. By corollary \ref{corollary f^* has closed graph} $gr ( f^* ) \subseteq X \times Y$ is a closed set. The definition of the closure gives $cl \big( gr ( f ) \big) \subseteq gr ( f^* )$.
        \item[$\supseteq :$]
        Let $( x_0 , y_0 ) \in gr ( f^* )$. By definition of the star multifunction, we have $y_0 \in Z_{x_0}$, where $Z_{x_0}$ is the set of extended values of $f$ at $x_0$. Using condition \ref{Proposition Multi-Split continuity with nets}\ref{definition sequential multi-split continuity a} of the characterization of multi-split continuity using nets, we can find a convergent net $\langle x_\lambda \rangle_{\lambda \in \Lambda} \longrightarrow x_0$ in $X$ with $\langle f ( x_\lambda ) \rangle_{\lambda \in \Lambda} \longrightarrow y_0$ in $Y$. In particular, we have the convergent net $\langle \big( x_\lambda , f ( x_\lambda ) \big) \rangle_{\lambda \in \Lambda} \longrightarrow ( x_0 , y_0 )$ in $X \times Y$. Note that $\big( x_\lambda , f ( x_\lambda ) \big) \in gr ( f )$ for all $\lambda \in \Lambda$. By the characterization of closed sets using nets in proposition \ref{proposition nets and closure}, the limit point $( x_0 , y_0 )$ therefore lies in the closure $cl \big( gr ( f ) \big)$ and we have $cl \big( gr ( f ) \big) \supseteq gr ( f^* )$.
    \end{description}
    This shows the second inclusion and thus concludes the proof.
\end{proof}

In the following remark, we present a method of constructing a multi-split continuous function $f_{weird} : [0 , 1 ] \to [ 0 , 1]$ with the property that $| f_{weird}^* ( \cdot ) |$ is locally unbounded at $0$. The proof of the multi-split continuity of this function is left to the reader.

\begin{remark}
    For every $n \in N$, choose $n+1$ distinct sequences \linebreak\mbox{$( x_k^{n , i} )_{k \in \mathbb{N} } \longrightarrow 1/n$} for $ i = 0 , ... \, , n$, where all $x_k^{n , i} \in B_d \big(1/ n(n+1), 1/n \big)$ are pairwise distinct for all possible choices of $k\in \mathbb{N}$ and $i = 0, ... \, , n$. This is possible, since $B_d \big(1/ n(n+1), 1/n \big) \subseteq \mathbb{R}$ is uncountable. Now define $f_{weird}$ as follows.
    \begin{align*}
        f_{weird} : [ 0 , 1 ] \to \mathbb{R}^2
        , \quad
        x \mapsto 
        \begin{cases}
            (x , i / n^2 )  \quad & \text{if } x = x_k^{n , i} \text{ for some } k , n , i \\
            (x , 0 ) & \text{else}
        \end{cases}
    \end{align*}
    Note that $x = x_k^{n , i}$ for at most one combination of the $k, n, i$. Thus, $f_{weird}$ is well-defined.
    
    Further note that $| f_{weird}^* ( x ) |$ is locally unbounded at 0, since we have \linebreak$f_{weird}^* ( 0 ) = \{ (0 , 0) \}$ and $f_{weird}^* ( 1 / n ) \, = \, \{ ( 1 / n , i / n^2 ) : i = 0 , ... \, , n \}$ for $n \in \mathbb{N}$. We therefore have $| f_{weird}^* ( 1 / n ) | = n + 1 \overset{n \to \infty}{\longrightarrow} \infty$.
\end{remark}

Recall the Kuratowski theorem from general topology. Given a Hausdorff space $Y$, the usual projection map $\pi_X : X \times Y \to X$ is closed for any space $X$ if and only if the space $Y$ is compact \cite[Thm.~3.1.16]{En89}. The following argument in part follows this proof.

\begin{proposition}
\label{proposition coordinate projection from graph of f^* is closed}
    Let $f : X \to Y$ be a globally multi-split continuous function, where $Y$ is a locally compact regular Hausdorff space, and let \linebreak$\pi_X : X \times Y \to X$ be the usual coordinate projection map. Then the restriction $\pi_X |_{gr ( f^* )} : gr ( f^* ) \to X$ is a closed map.
\end{proposition}

\begin{proof}
    Let $C \subseteq gr ( f^* )$ be an arbitrary closed subset. We will show that \linebreak$\widetilde{C} := \pi_X |_{gr ( f^* )} ( C ) \subseteq X$ is closed. Let $x_0 \in X \, \backslash \, \widetilde{C}$. We will show that there exists an open neighborhood $U \subseteq X$ of $x_0$ with $U \cap \widetilde{C} = \varnothing$. This is equivalent to having $( U \times Y ) \cap C = \varnothing$.
    
    The set of extended values $Z_{x_0} \subseteq Y$ of $f$ at $x_0$ is finite by definition. The space $Y$ is locally compact, and therefore there exists a compact neighborhood $K \subseteq Y$ of $Z_{x_0}$.

    By definition $C$ is a closed subset of $gr ( f^* )$. Corollary \ref{corollary f^* has closed graph} states that $gr ( f^* )$ is closed in $X \times Y$ and $C$ is therefore also closed as a subset of $X \times Y$. We have $( \{ x_0 \} \times K ) \cap C = \varnothing$. For every $y \in K$, the point $( x_0 , y ) \in C^c$ and there thus exist open neighborhoods $U_y \subseteq X$ of $x_0$ and $V_y \subseteq Y$ of $y$ in $Y$ such that $( U_y \times V_y ) \cap C = \varnothing$. Then $( V_y )_{ y \in K }$ is an open covering of the compact set $K$, and there exists a finite subcovering $V_{y_1} , ... \, , V_{y_N}$. Let $U := \bigcap_{i=1}^{N} U_{y_i} \subseteq X$. Then $U$ is an open neighborhood of $x_0$ in $X$ with $( U \times K ) \cap C = ( U \times Y ) \cap C = \varnothing$. This concludes the proof.
\end{proof}

Recall yet another fact from general topology. Let $f : X \to Y$ be a continuous map, where $X$ is compact and $Y$ is Hausdorff. Then $f$ is a closed map \cite[Thm.~3.1.12]{En89}.

\begin{proposition}
    Let $X$ be a compact Hausdorff space, let $Y$ be regular Hausdorff, and let $f : X \to Y$ be globally multi-split continuous. Then the restriction of the usual coordinate projection map $\pi_X |_{gr ( f^* )} : gr ( f^* ) \to X$ is closed.
\end{proposition}

\begin{proof}
    The star multifunction $f^*$ is upper semicontinuous by theorem \ref{theorem multi-split extension is usc}. We see that $f^* ( X ) \subseteq Y$ is compact by corollary \ref{corollary f^* sends cpt sets to cpt sets}. Corollary \ref{corollary f^* has closed graph} implies that $gr ( f^* ) \subseteq X \times Y$ is closed. In particular, $gr ( f^* )$ is a closed subset of the compact space $X \times f^* ( X )$, and therefore $gr ( f^* )$ is a compact space itself. The map $\pi_X |_{gr ( f^* )} : gr ( f^* ) \to X$ then is a continuous map with a compact domain and a Hausdorff codomain. Such maps are closed by the statement quoted above.
\end{proof}

\section{Pre-Multi-Split Multifunctions}
\label{Section Pre-multi-split Functions}
\noindent
In this section, we will study multifunctions whose selections are multi-split continuous, referring to them as pre-multi-split. Theorem \ref{theorem finite usc functions are pre-multi-split} states that all non-empty, finite, upper semicontinuous multifunctions are pre-multi-split. Proposition \ref{proposition f^* is pre-multi-split} extends this result to star multifunctions with regular Hausdorff domains. Corollary \ref{corollary union of fin. many multi-split functions is pre-multi-split} establishes that finite unions of multi-split continuous functions are pre-multi-split. These results provide a systematic way to construct multi-split continuous functions without having to verify the property at each point individually.

\begin{definition}
\label{Definition pre-multi-split}
    A multifunction $F : X \rightrightarrows Y$ is called \textit{pre-multi-split at $p \in X$} if every selection $f : X \to Y$ of $F$ is multi-split continuous at $p \in X$. We call $F$ \textit{pre-multi-split} if every selection $f$ of $F$ is globally multi-split continuous.
\end{definition}

Pre-multi-split multifunctions have the following connection with subcontinuous multifunctions.

\begin{proposition}
    Let $F : X \rightrightarrows Y$ be a multifunction that is subcontinuous at $p \in X$. We define $\widetilde{Z}_{p} \subseteq Y$ to be the following set.
    \begin{equation*}
        \widetilde{Z}_{p} := 
        \left\{
        y \in Y \;\middle|\;
        \parbox{0.57\textwidth}{
            There exists a convergent net $\langle x_\lambda \rangle_{\lambda \in \Lambda} \longrightarrow p$ in $X$ such that $y$ is a cluster point of the net $\langle y_\lambda \rangle_{\lambda \in \Lambda}$ with $y_\lambda \in F (x_\lambda)$ for all $\lambda \in \Lambda$ }
        \right\}
    \end{equation*}
     If $\widetilde{Z}_{p}$ is finite, then $F$ is pre-multi-split at $p$.
\end{proposition}

\begin{proof}
    Let $f : X \to Y$ be a selection of $F$. We define the set $Z_p$ to be the following subset of $\widetilde{Z}_p$ in $Y$.
    \begin{equation*}
        Z_{p} := 
        \left\{
        y \in Y \;\middle|\;
        \parbox{0.65\textwidth}{
            There exists a convergent net $\langle x_\lambda \rangle_{\lambda \in \Lambda} \longrightarrow p$ in $X$ such that $y$ is a cluster point of the net $\langle f ( x_\lambda ) \rangle_{\lambda \in \Lambda}$}
        \right\}
    \end{equation*}
    We clearly have $f ( p ) \in Z_p$ and thus $Z_p$ is non-empty. Further, $Z_p$ is a subset of a finite set and therefore itself finite.
    
    We now verify that $Z_p$ fulfills conditions \ref{Proposition Multi-Split continuity with nets}\ref{definition sequential multi-split continuity a} and \ref{Proposition Multi-Split continuity with nets}\ref{definition sequential multi-split continuity b} of the characterization of multi-split continuity using nets. The cluster points of a net coincide with the possible limit points of its subnets. Condition \ref{Proposition Multi-Split continuity with nets}\ref{definition sequential multi-split continuity a} therefore holds trivially for any $y_i \in Z_p$. By assumption, the multifunction $F$ is subcontinuous at $p$. In particular, for any convergent net $\langle x_\lambda \rangle_{\lambda \in \Lambda} \longrightarrow p$ in $X$, there exists a cluster point of the net $\langle f ( x_\lambda ) \rangle_{\lambda \in \Lambda }$ which clearly lies in $Z_p$. Condition \ref{Proposition Multi-Split continuity with nets}\ref{definition sequential multi-split continuity b} therefore also follows directly.
\end{proof}

\begin{theorem}
\label{theorem finite usc functions are pre-multi-split}
    Let $F : X \rightrightarrows Y$ be a multifunction that is non-empty, finite, and upper semicontinuous (at $p \in X$). Then $F$ is pre-multi-split (at $p \in X$).
\end{theorem}

\begin{proof}
    Let $g : X \to Y$ be a selection of $F$. We verify that $g$ fulfills condition \ref{proposition equivalent definition only b 2} in proposition \ref{proposition equivalent definition only with condition b} using the set $Z := F ( p )$. The multifunction $F$ is non-empty and finite, and thus so is the set $Z$. Let $V \subseteq Y$ be an open neighborhood of $Z$. The multifunction $F$ is upper semicontinuous at $p$. Thus, there exists a neighborhood $U \subseteq X$ of $p$ with $g ( U ) \subseteq F ( U ) \subseteq V$. Using proposition \ref{proposition equivalent definition only with condition b}, we may conclude that $g$ is multi-split continuous at $p$. Therefore, $F$ is pre-multi-split at $p \in X$.
\end{proof}

\begin{proposition}
\label{proposition f^* is pre-multi-split}
    Let $Y$ be a regular Hausdorff space and let $f \in mSplit ( X , Y )$ be a globally multi-split continuous function. Then $f^* : X \rightrightarrows Y$ is pre-multi-split.
\end{proposition}

\begin{proof}
    The star multifunction $f^* : X \rightrightarrows Y$ is globally upper semicontinuous by theorem \ref{theorem multi-split extension is usc}, since we assumed $Y$ to be a regular Hausdorff space. Further, by definition, $f^* ( p )$ is finite for all $p \in X$. We may therefore apply theorem \ref{theorem finite usc functions are pre-multi-split} to conclude this proof.
\end{proof}

Note that not every pre-multi-split function $F : X \rightrightarrows Y$ is of the form $f^*$ for some $f \in mSplit ( X , Y )$. Consider $F : \{ p \} \rightrightarrows \mathbb{R}, \, \allowbreak p \mapsto \mathbb{R}$. The following proposition at least gives a sufficient condition.

\begin{proposition}
\label{proposition minimal usco gives G = f^*}
    Let $Y$ be regular Hausdorff and let $G : X \rightrightarrows Y$ be a finite minimal usco map. Then every selection $f : X \to Y$ of $G$ satisfies $f^* \equiv G$.
\end{proposition}

\begin{proof}
    Let $f : X \to Y$ be a selection of $F$. By proposition \ref{proposition closed graph theorem for usc multifunctions}, the graph $gr ( F ) \subseteq X \times Y$ of $F$ is closed. In particular, we have $cl \big( gr ( f ) \big) \subseteq gr ( F )$. Proposition \ref{proposition graph of f^* is closure of graph of f} states $cl \big( gr ( f ) \big) = gr ( f^* )$. We thus have $gr ( f^* ) \subseteq gr ( F )$. Further, $f^*$ is an usco map by corollary \ref{corollary f^* is usco map}. The multifunction $F$ was assumed to be minimal usco, and we therefore have $F \equiv f^*$ as stated.
\end{proof}

Note that the star multifunction of the characteristic function $\mathds{1}_\mathbb{Q} : \mathbb{R} \to \mathbb{R}$ is not minimal usco. The condition that the multifunction $G$ is minimal usco can therefore not be necessary for the existence of a function $f$ with $f^* \equiv G$.

The following proposition generalizes a statement by the second author \cite[Prop.~3.13]{Gh24}, the proof presented here is an adaptation.

\begin{proposition}
    Let $f : X \to Y$ be a function and $Y$ be a Hausdorff space. Then $f$ is continuous if and only if $f$ is globally multi-split continuous and has a closed graph.
\end{proposition}

\begin{proof}
    We verify two implications. Let the space $Y$ be Hausdorff and let $f : X \to Y$ be a single-valued function. We identify $f$ with the multifunction $f : X \rightrightarrows Y, \, x \mapsto \{ f ( x ) \}$.

    First, let $f : X \to Y$ be continuous. By proposition \ref{proposition subcontinuous and closed graph iff continuous}, we see that $f$ is subcontinuous and has a closed graph. Proposition \ref{proposition connection of upper semicontinuity and subcontinuity} then implies that the multifunction $f$ is upper semicontinuous, and by theorem \ref{theorem finite usc functions are pre-multi-split}, $f$ is pre-multi-split. The single-valued function $f$ therefore is multi-split continuous as stated.

    Now let $f$ be globally multi-split continuous. By proposition \ref{proposition graph of f^* is closure of graph of f}, we have that $gr ( f^* ) = cl \big( gr ( f ) \big) = gr ( f )$ and thus $f^* \equiv f$. Theorem \ref{theorem multi-split extension is usc} then implies that $f^* \equiv f$ is upper semicontinuous. We apply proposition \ref{proposition connection of upper semicontinuity and subcontinuity} to see that $f$ is subcontinuous and $f$ is continuous by proposition \ref{proposition subcontinuous and closed graph iff continuous}.
\end{proof}

We formulate yet another equivalent definition of global multi-split continuity of $f : X \to Y$, where $Y$ is now assumed to be a compact regular Hausdorff space.

\begin{proposition}
\label{Proposition equivalent definition using closure of the graph}
    Let $Y$ be a compact regular Hausdorff space and $f : X \to Y$ be a function. Let $Z := cl \big( gr ( f ) \big) \subseteq X \times Y$ denote the closure of the graph of $f$ in $X \times Y$. Then the following two conditions are equivalent.
    \begin{enumerate}[label=\arabic*), nosep]
        \item
        $f : X \to Y$ is globally multi-split continuous

        \item
        The standard projection $\pi : X \times Y \to X$ restricted to $Z$ is finite-to-one.
    \end{enumerate}
\end{proposition}

\begin{proof}
    We show two implications. First let $f : X \to Y$ be globally multi-split continuous and $p \in X$ be arbitrary. By proposition \ref{proposition graph of f^* is closure of graph of f}, we have that $Z = gr ( f^* )$. Let $Z_p \subseteq Y$ be the set of extended values of $f$ at $p \in X$. We then have
    $$( \pi |_Z )^{-1} \big( \{ p \} \big) = \{ p \} \times Z_{p}$$
    The set $Z_p \subseteq Y$ is finite by definition. Thus, $\pi |_Z$ is finite-to-one.
    
    Now assume that the map $\pi |_Z : cl \big( gr ( f ) \big) \to X$ is finite-to-one. We define $F : X \rightrightarrows Y$ to be the following multifunction.
    $$F : X \rightrightarrows Y, \quad p \mapsto \{ y \in Y : ( p , y ) \in cl \big( gr ( f ) \big) \}$$
    $F$ is finite and non-empty. Further, its graph $gr ( F ) = Z$ is closed by assumption. Applying proposition \ref{proposition closed graph theorem for usc multifunctions}, we see that $F$ is a finite upper semicontinuous multifunction and thus pre-multi-split by theorem \ref{theorem finite usc functions are pre-multi-split}. The function $f$ is multi-split continuous as a selection of $F$, and we may conclude the second implication.
\end{proof}

Proposition \ref{Proposition equivalent definition using closure of the graph} above might fail if we allow the space $Y$ to be non-{\linebreak}sequentially compact. This is demonstrated by the following example.

\begin{example}
    We define the following function $f : [ 0 , 1 ] \to ( 0 , 1 ]$.
    $$f : [ 0 , 1 ] \to ( 0 , 1 ] , \quad
    x \mapsto 
    \begin{cases}
        1/n \quad &\text{if } x = 1 / n \text{ for some } n \in \mathbb{N} \\
        1 \quad &\text{else}
    \end{cases}$$
    The closure of the graph of $f$, which we denote by $Z$, is then given by the following expression.
    $$Z \, := \, cl \big( gr ( f ) \big) \, = \, \big( [0 , 1] \times \{ 1 \} \big) \cup \{ ( 1/n , 1/n ) : n \in \mathbb{N} \}$$
    Clearly, the projection $\pi_X : [ 0 , 1 ]\times ( 0 , 1 ] \to [ 0 , 1 ]$ onto the first coordinate is at most two-to-one when restricted to $Z$. However, we can construct the convergent net $\langle 1 / n \rangle_{n \in \mathbb{N} } \longrightarrow 0$ in $[ 0 , 1 ]$ where the net $\langle f ( 1 / n ) \rangle_{n \in \mathbb{N} } = \langle 1/n \rangle_{n \in \mathbb{N}}$ does not have any cluster points in $( 0 , 1 ]$. This contradicts condition \ref{Proposition Multi-Split continuity with nets}\ref{definition sequential multi-split continuity b} of the characterization of multi-split continuity using nets. Therefore $f$ is not multi-split continuous at $0$.
    
\end{example}

\begin{lemma}
\label{Lemma unions of multi-split functions are pre-multi-split}
    Let $f , g \in mSplit ( X , Y )$ be multi-split continuous functions. Then the multifunction \mbox{$f \cup g : X \rightrightarrows Y$} is pre-multi-split.
\end{lemma}

\begin{proof}
    Let $h : X \to Y$ be a selection of the multifunction $f \cup g : X \rightrightarrows Y$ and let $p \in X$ be arbitrary. We show multi-split continuity of $h$ at $p$ using proposition \ref{proposition equivalent definition only with condition b}. Let $Z_p^f \subseteq Y$ and $Z_p^g \subseteq Y$ be sets of extended values of $f$ and $g$ at $p$. We define $Z := Z_p^f \cup Z_p^g$. Note that $Z \subseteq Y$ is non-empty and finite. Now let $V \subseteq Y$ be a neighborhood of $Z$. In particular, $V$ is a neighborhood of $Z_p^f$ in $Y$, and by condition \ref{Definition topological multi-split function}\ref{definition topological multi-split continuity b} of the multi-split continuity of $f$ at $p$, there exists a neighborhood $U_f \subseteq X$ of $p$ with $f ( U_f ) \subseteq V$. Analogously, we can find a neighborhood $U_g \subseteq X$ of $p$ with $g ( U_g ) \subseteq V$. We set $U := U_f \cap U_g$ and then have the following inclusions.
    $$h ( U ) \, \subseteq \, f ( U ) \cup g ( U ) \, \subseteq \, f ( U_f ) \cup g ( U_g ) \, \subseteq \, V$$
    Thus, the function $h$ satisfies condition \ref{proposition equivalent definition only b 2} in proposition \ref{proposition equivalent definition only with condition b}, and $h$ is multi-split continuous at $p \in X$.
\end{proof}

\begin{proposition}
\label{proposition properties of pre-multi-split functions}
    Let $F, F_1 , ... \, , F_n : X \rightrightarrows Y$ be pre-multi-split multifunctions. Then the following multifunctions are also pre-multi-split.
    \begin{enumerate}[label=\arabic*), left=0pt, nosep]
        \item
        \label{proposition properties of pre-multi-split functions property 1}
        Any non-empty submultifunction $G : X \rightrightarrows Y$ of $F$.

        \item
        \label{proposition properties of pre-multi-split functions property 2}
        Any finite union $G := F_1 \cup \cdots \cup F_n : X \rightrightarrows Y$ of pre-multi-split functions.
    \end{enumerate}
\end{proposition}

\begin{proof}
    Any selection of a submultifunction is also a selection of the original multifunction. The statement \ref{proposition properties of pre-multi-split functions property 1} therefore holds trivially. The proof of statement \ref{proposition properties of pre-multi-split functions property 2} requires some more work.

    By induction, it suffices to show the statement in the case $n = 2$. Thus let $F_1 , F_2 : X \rightrightarrows Y$ be pre-multi-split multifunctions and let $f : X \to Y$ be a selection of $F_1 \cup F_2$. We show that $f$ is multi-split continuous. We define the two functions $g , h : X \to Y$ using the following procedure.
    $$
    g : X \to Y, \quad x \mapsto
    \begin{cases}
        f ( x ) &\text{if } f(x) \in F_1 (x) \\
        y_x \in F_1 ( x ) &\text{else}
    \end{cases}
    $$
    $$
    h : X \to Y, \quad x \mapsto
    \begin{cases}
        f ( x ) &\text{if } f(x) \in F_2 (x) \\
        y_x \in F_2 ( x ) &\text{else}
    \end{cases}
    $$
    Here the values $y_x \in F_i ( x )$ are chosen arbitrarily from the non-empty sets $F_i ( x )$. Note that $g$ is a selection of $F_1$ and $h$ is a selection of $F_2$. By the definition of a pre-multi-split multifunction, both $g$ and $h$ are therefore globally multi-split continuous. Further note that $f$ is a selection of $g \cup h$, since ${f ( x ) \in F_1 ( x ) \cup F_2 ( x )}$ for all $x \in X$. We may therefore apply lemma \ref{Lemma unions of multi-split functions are pre-multi-split} and conclude the proof.
\end{proof}

Note that arbitrary intersections $\bigcap_{i \in I} F_i : X \rightrightarrows Y$ of multifunctions are submultifunctions of any of the $F_i$. The statement \ref{proposition properties of pre-multi-split functions property 2} therefore is also true for arbitrary non-empty intersections of pre-multi-split multifunctions.

We can now describe a method of constructing a wide variety of multi-split continuous functions from a given set of multi-split continuous functions in the following corollary of proposition \ref{proposition properties of pre-multi-split functions}.

\begin{corollary}
\label{corollary union of fin. many multi-split functions is pre-multi-split}
    Given a finite collection of multi-split continuous functions \linebreak$f_1 , ... \, , f_n : X \to Y$ we can construct a new multi-split continuous function $g : X \to Y$ by taking a selection of the multifunction $\bigcup_{i=1}^n f_i : X \rightrightarrows Y$.
\end{corollary}

\begin{remark}
\label{remark compositions of pre-multi-split functions are not pre-multi-split}
    We define the two multifunctions $F : \mathbb{R} \rightrightarrows \{ * \} , \, x \mapsto \{ * \}$ and $G : \{ * \} \rightrightarrows \mathbb{R}, \, * \mapsto \mathbb{R}$, where $\{ * \}$ is a single-point space. Note that any selection of $F$ is constant and thus continuous. Also, any selection of $G$ is a function with a discrete domain and therefore continuous. Both multifunctions $F$ and $G$ are therefore pre-multi-split. By definition, their composition is $G \circ F : \mathbb{R} \rightrightarrows \mathbb{R} , \, x \mapsto \mathbb{R}$, which clearly is not pre-multi-split. We have thus demonstrated that compositions of pre-multi-split functions need not be pre-multi-split.
\end{remark}

We now turn towards inverse-image multifunctions and study when they are pre-multi-split. This will be particularly relevant in section \ref{Section Cuts and Re-glues}.

\begin{proposition}
\label{proposition inverse image multifunction is pre-multi-split}
    Let $f : X \to Y$ be a closed finite-to-one surjection. Then the inverse-image multifunction $f^{-1} : Y \rightrightarrows X$ is pre-multi-split.
\end{proposition}

\begin{proof}
    The function $f$ is a finite-to-one surjection. We immediately get that its inverse image multifunction $f^{-1} : Y \rightrightarrows X$ is non-empty and finite. Using proposition \ref{proposition inverse image multifunction in u.s.c. iff f is closed}, we see that $f^{-1}$ is an upper semicontinuous multifunction. We may therefore apply theorem \ref{theorem finite usc functions are pre-multi-split} to conclude.
\end{proof}

One follow-up question is: When is the inverse image multifunction of a continuous quotient projection with finite equivalence classes pre-multi-split? In the following propositions and examples, we present a sufficient condition. 

\begin{proposition}
\label{proposition inverse image multifunction of projection from compact to Hausdorff is pre-multisplit}
    Let $X$ be a compact space and let $\sim$ be an equivalence relation on $X$ whose equivalence classes are finite such that the space $X /_\sim$ is Hausdorff. Let $p : X \to X /_\sim$ be the usual quotient map. Then the inverse image multifunction $p^{-1} : X /_\sim \rightrightarrows X$ is pre-multi-split.
\end{proposition}

\begin{proof}
    The map $p : X \to X /_\sim$ is continuous with a compact domain and a Hausdorff codomain and thus closed \cite[Thm.~3.1.12]{En89}. We may apply proposition \ref{proposition inverse image multifunction is pre-multi-split} and conclude.
\end{proof}

We now present a counterexample to proposition \ref{proposition inverse image multifunction of projection from compact to Hausdorff is pre-multisplit}, where the space $X$ is not sequentially compact.

\begin{example}
    Consider the space $X := \mathbb{R}_{\geq 0 }$ and define the equivalence relation $\sim$ on $\mathbb{R}_{\geq 0}$.
    $$x \sim y 
    \quad : \Longleftrightarrow \quad
    x = y \quad \text{or} \quad \{ x , y \} = \{ n , 1/n \} \text{ for some } n \in \mathbb{N}$$
    Let $p_\sim : X \to X /_\sim$ be the usual quotient map. The equivalence classes of $\sim$ contain at most two elements and therefore are finite. It is left to the reader to verify that $X /_\sim$ is in fact Hausdorff. We set $f$ to be the following selection of $p_\sim^{-1} : X /_\sim \rightrightarrows X$.
    $$f : X /_\sim \to X , \quad
    [ x ]_\sim \mapsto 
    \begin{cases}
        n \quad &\text{if} \quad n \in [ x ]_\sim \text{ for some } n \in \mathbb{N}\\
        x \quad &\text{else}
    \end{cases}$$
    We then have the convergent net $\langle [ 1/n ]_\sim \rangle_{n \in \mathbb{N} } \longrightarrow [0]_\sim$ in $X /_\sim$ where the net $\langle f ( [ 1 / n ]_\sim ) \rangle_{n \in \mathbb{N} } = \langle n \rangle_{n \in \mathbb{N} }$ has no cluster points in $X$. This contradicts condition \ref{Proposition Multi-Split continuity with nets}\ref{definition sequential multi-split continuity b} of the characterization of multi-split continuity using nets. Therefore $f$ is not multi-split continuous and $p_\sim^{-1}$ is not pre-multi-split.
\end{example}

We further outline a counterexample to proposition \ref{proposition inverse image multifunction of projection from compact to Hausdorff is pre-multisplit}, where we drop the assumption that $X /_\sim$ is Hausdorff.

\begin{example}
    We set $X \subseteq \mathbb{R}^2$ to be the following space.
    $$X := \{ ( 0 , 0 ) \} \, \cup \, \bigcup_{n=1}^\infty \, [ 0 , 1 / n ] \times \{ 1 / n \} \subseteq \mathbb{R}^2$$
    Clearly, $X$ is a bounded subset. Further, one can easily see that $X$ is closed. Using the Heine--Borel theorem we may conclude that $X$ is compact. We now define the equivalence relation $\sim$ on $X$.
    $$(x_0 , y_0 ) \sim (x_1 , y_1 )
    \quad : \Longleftrightarrow \quad
    ( x_0, y_0 ) = ( x_1 , y_1 ) \quad \text{or} \quad x_0 = x_1 \neq 0
    $$
    We have the following expression for the equivalence class of any element $( x_0 , y_0 ) \in X$. Note that all equivalence classes are finite.
    $$[ ( x_0 , y_0 ) ]_\sim = \{ ( x_0 , 1 / k ) : 1 / k \geq y_0 \}  \text{ if } x_0 \neq 0
    \quad \text{and} \quad
    [ ( 0 , y_0 ) ]_\sim = \{ ( 0 , y_0 ) \}
    $$
    Choose pairwise disjoint sequences of distinct values $( x_k^n )_k \longrightarrow 0$ for every \linebreak$n \in \mathbb{N}$ with $x_k^n \in ( 0 , 1/n )$. That is, let $x_{k_1}^{n_1} = x_{k_2}^{n_2}$ if and only if $k_1 = k_2$ and $n_1 = n_2$ for all $n_i , k_i \in \mathbb{N}$. This is possible, since every interval $( 0 , 1/n )$ is uncountable. Now define the following right inverse $f : X /_\sim \to X$ of the quotient map $p_\sim : X \to X_\sim$.
    $$f : X /_\sim \to X , \quad
    [ ( x , y ) ]_\sim \mapsto
    \begin{cases}
        ( 0 , y ) \quad &\text{if } x = 0 \\
        ( x , 1 / n ) \quad &\text{if } y = x^n_k \text{ for some } k , n \in \mathbb{N} \\
        ( x , 1 ) \quad &\text{else}
    \end{cases}$$
    The net $\langle [ ( x^n_k , 1 ) ]_\sim \rangle_{ ( n , k ) \in \mathbb{N}^2 }$ then converges to any of the $[ (0 , 1/n ) ]_\sim$, where $n \in \mathbb{N}$. The upward directed preorder $\preccurlyeq$ on $\mathbb{N}^2$ is defined as follows.
    $$( n_1 , m_1 ) \preccurlyeq ( n_2 , m_2 ) 
    \quad : \Longleftrightarrow \quad
    n_1 \leq n_2 \text{ and } m_1 \leq m_2$$
    We then have $\langle f \big( [ ( x_k^n , 1 ) ]_\sim \big) \rangle_{(k , n ) \in \mathbb{N}^2} = \langle ( x_k^n , 1/n ) \rangle_{(k , n ) \in \mathbb{N}^2}$. Clearly, every \linebreak$(0 , 1/n ) \in X$ is a cluster point of this net. By proposition \ref{proposition net of function values can only have fin. many cluster points} the function $f$ cannot be split continuous at any of the $[ ( 0 , 1/n ) ]_\sim$ since for any $n \in \mathbb{N}$, the set $f^* \big( [ ( 0 , 1/n ) ]_\sim \big)$ would have to contain the infinite set $\{ ( 0 , 1/m ) : m \in \mathbb{N} \}$, a contradiction.
\end{example}

\section{Split Homeomorphisms}
\label{Section Split Homeomorphisms}
\noindent
We now extend the concept of homeomorphisms by introducing split homeomorphisms. The stability of multi-split continuous functions under composition allows us to establish split homeomorphism of spaces as an equivalence relation on topological spaces, as formalized in Proposition \ref{proposition split homeomorphism is an equivalence relation}. Corollary \ref{corollary split homeomorphic spaces preserve compactness} states that split homeomorphisms preserve compactness, demonstrating the non-triviality of this notion.

\begin{definition}
\label{Definition split Homeomorphisms}
    Let $f \in mSplit ( X , Y )$ be a multi-split continuous bijection. We call $f$ a \textit{split homeomorphism} if its inverse $f^{-1} : Y \to X$ is multi-split continuous as well. We call the spaces $X$ and $Y$ \textit{split homeomorphic} and write $X  \Bumpeq Y$. We further set $mHomeo (X , Y )$ to be the set of all split homeomorphisms $f : X \to Y$.
\end{definition}

\begin{remark}
\label{Remark homeomorphic spaces are split homeomorphic}
    Continuous functions are multi-split continuous by remark \ref{Remark multi-split versus continuous and split-continuous}. Homeomorphisms, therefore, also are split homeomorphisms.
\end{remark}

\begin{proposition}
\label{proposition split homeomorphism is an equivalence relation}
    Split homeomorphism of spaces is an equivalence relation. 
\end{proposition}

\begin{proof}
    Let $X$, $Y$, and $Z$ be arbitrary topological spaces. As mentioned in remark \ref{Remark homeomorphic spaces are split homeomorphic}, the homeomorphism $id_X : X \to X$ also is a split homeomorphism. We thus have $X \Bumpeq X$ and $\Bumpeq$ is reflexive. Symmetry of $\Bumpeq$ is a direct consequence of the definition, since $( f^{-1} )^{-1} \equiv f$. Finally, the transitivity of $\Bumpeq$ directly follows from theorem \ref{Theorem compositions of multi-split functions} above, stating that compositions of multi-split continuous functions are multi-split continuous. Thus $\Bumpeq$ is an equivalence relation.
\end{proof}

\begin{example}
    Let $X$ be a finite set and let $\tau_1 , \tau_2$ be topologies on $X$. We define the following function.
    $$\varphi : ( X , \tau_1 ) \to ( X , \tau_2 ) , \quad x \mapsto x$$
    By applying proposition \ref{proposition equivalent definition only with condition b} to $\varphi$ and to $\varphi^{-1}$ with $Z := X$, we immediately see that $\varphi$ is a split homeomorphism. In particular, we have the following split homeomorphic spaces.
    $$( X , \tau_1 ) \Bumpeq ( X , \tau_{triv} ) \Bumpeq ( X , \tau_{disc} )$$
    This demonstrates that neither connectedness nor axioms of separation of topological spaces are preserved under split homeomorphisms. Compared to homeomorphisms of spaces, split homeomorphisms therefore define a truly coarser equivalence relation on the class of all topological spaces.
\end{example}

\begin{proposition}
    Let $X$ and $Y$ be regular Hausdorff spaces, and let \linebreak$f : X \to Y$ be a split homeomorphism. Then for every $y_0 \in Y$ we have the following expression of $(f^{-1} )^* ( y_0 )$.
    $$( f^{-1} )^* ( y_0 ) = \{ x \in X : ( x , y_0 ) \in gr ( f^* ) \}$$
\end{proposition}

\begin{proof}
    Let $y_0 \in Y$ be an arbitrary point. Applying proposition \ref{proposition graph of f^* is closure of graph of f} to the multi-split continuous functions $f : X \to Y$ and $f^{-1} : Y \to X$ we obtain the following expressions for $gr ( f^* )$ and $gr \big( ( f^{-1} )^* \big)$.
    $$gr ( f^* ) = cl \big( gr ( f ) \big) \subseteq X \times Y 
    \quad \text{and} \quad 
    gr \big( ( f^{-1} )^* \big) = cl \big( gr ( f^{-1} ) \big) \subseteq Y \times X$$
    Let $\cdot^T : X \times Y \to Y \times X, \, ( x , y ) \mapsto ( y , x)$ be the usual reflection map. Note that the reflection map $\cdot^T$ commutes with the closure operation on subsets in the product topology. Further, note that we have $gr ( f )^T = gr ( f^{-1} )$. We can therefore write the following equivalences.
    \begin{equation*}
        gr \big( (f^{-1} )^* \big) = cl \big( gr ( f^{-1} ) \big) = cl \big( gr ( f )^T \big) = cl \big( gr ( f ) \big)^T = gr ( f^* )^T
    \end{equation*}
    For any $x \in X$ we now have $x \in ( f^{-1} )^* ( y_0 )$ if and only if $( y_0 , x ) \in gr \big( (f^{-1} )^* \big)$. By the calculation above, that is true if and only if $( y_0 , x ) \in gr ( f^* )^T$ or $( x , y_0 ) \in gr ( f^* )$ as claimed in the statement above.
\end{proof}

We now observe that compactness of topological spaces is preserved under split homeomorphism, which makes this notion nontrivial.

\begin{proposition}
\label{Proposition split homeomorphisms and compact spaces}
    Let $f : X \to Y$ be a multi-split continuous surjection. If $X$ is compact, then so is $Y$.
\end{proposition}

\begin{proof}
    We use the characterization of compact spaces using nets as seen in proposition \ref{proposition nets and compact spaces}. Let $\langle y_\lambda \rangle_{\lambda \in \Lambda}$ be an arbitrary net in $Y$. The function $f$ is surjective. For every $\lambda \in \Lambda$ we can therefore choose some $x_\lambda \in X$ with $f ( x_\lambda ) = y_\lambda$. The space $X$ is compact, and therefore the net $\langle x_\lambda \rangle_{\lambda \in \Lambda}$ has some cluster point $p \in X$. By proposition \ref{proposition cluster points of a net are exactly limits of convergent subnets}, there exists a convergent subnet $\langle x_{\phi ( \theta )} \rangle_{\theta \in \Theta} \longrightarrow p$. The function $f$ is multi-split continuous at $p$ by assumption. Condition \ref{Proposition Multi-Split continuity with nets}\ref{definition sequential multi-split continuity b} of the characterization of multi-split continuous functions using nets then states that the net $\langle f ( x_{\phi ( \theta )} ) \rangle_{\theta \in \Theta} = \langle y_{\phi ( \theta )} \rangle_{\theta \in \Theta}$ has a cluster point $y_0 \in Y$. Applying proposition \ref{proposition cluster points of a net are exactly limits of convergent subnets} once more gives that $y_0$ also is a cluster point of the net $\langle y_\lambda \rangle_{\lambda \in \Lambda}$. We have thus found a cluster point of an arbitrary net in $Y$ and may conclude the proof using proposition \ref{proposition nets and compact spaces}.
\end{proof}

As a direct consequence of proposition \ref{Proposition split homeomorphisms and compact spaces}, we can now state the following corollary.

\begin{corollary}
\label{corollary split homeomorphic spaces preserve compactness}
    Let $X \Bumpeq Y$ be two split homeomorphic spaces. If $X$ is compact, then so is $Y$.
\end{corollary}

Recall the homeomorphism theorem from topology. Let $f : X \to Y$ be a continuous bijection, where $X$ is compact and $Y$ is Hausdorff. Then $f$ is a homeomorphism \cite[Thm.~3.1.13]{En89}. We now state a weaker version of this theorem for split homeomorphism.

\begin{proposition}
    Let $X$ be a compact space and $Y$ be a Hausdorff space. Let $f : X \to Y$ be a multi-split continuous bijection. If there exist finitely many continuous bijections $f_1 , ... \, , f_n : X \to Y$ such that $f ( p ) \in \bigcup_{i=1}^n f_i ( p )$ for all $p \in X$, then $f$ is a split homeomorphism.
\end{proposition}

\begin{proof}
    By applying the homeomorphism theorem to $f_i$ for each $i = 1 , ... \, , n$, we see that $f_i : X \to Y$ is a homeomorphism. In particular, $f_i^{-1} : Y \to X$ is continuous for every $i = 1 , ... \, , n$, and the multifunction $G : = \bigcup_{i=1}^n f_i^{-1}$ is pre-multi-split. The function $f^{-1}$ is a selection $G$ and therefore is multi-split continuous as stated.
\end{proof}

Despite a considerable effort, the authors were unable to prove or disprove the following conjecture.

\begin{conjecture}
\label{theorem split homeomorphism theorem}
    Given a compact space $X$ and a Hausdorff space $Y$, every multi-split continuous bijection $f : X \to Y$ is a split homeomorphism. That is, the inverse map $f^{-1} : Y \to X$ is multi-split continuous.
\end{conjecture}

\section{Cuts and with Subsequent Re-glues}
\label{Section Cuts and Re-glues}
\noindent
In this final section, we will explain the geometric intuition behind split homeomorphic spaces. We will first define a rigorous mathematical notion of two compact regular Hausdorff spaces being equivalent with cuts and subsequent re-glues. In propositions \ref{proposition equivalence with cuts and re-glues implies split homeomorphic} and \ref{proposition split homeomorphic implies equivalence with cuts and subsequent re-glues}, we will then show that two compact regular Hausdorff spaces are equivalent with cuts and subsequent re-glues if and only if they are split homeomorphic.

\begin{definition}
\label{definition equivalence with cuts and reglues}
    We say that two compact regular Hausdorff spaces $X$ and $Y$ are \textit{equivalent with cuts and subsequent re-glues} or \textit{$X$ can be deformed into $Y$ with cuts and subsequent re-glues} if there exists a compact space $Z$ together with the following data.
    \begin{enumerate}[label=\arabic*), nosep]
        \item 
        Two continuous finite-to-one surjections $p_X : Z \to X$ and $p_Y : Z \to Y$.
        \item
        A right inverse $p_X^{-1} : X \to Z$ of $p_X$ such that $f := p_Y \circ p_X^{-1}$ is a bijection.
    \end{enumerate}
    We write $X \bumpeq Y$. This situation can be described in the following diagram.
    \begin{center}
        \begin{tikzcd}
            &
                Z
                \ar[dl , "p_X"]
                \ar[dr , "p_Y"]
                \\
            
            X
            \ar[ur , "p_X^{-1}", bend left = 30]
            \ar[rr, "f"']
            &
                &
                    Y
        \end{tikzcd}
    \end{center}
\end{definition}

Recall the following fact from general topology. Let $f : X \to Y$ be a continuous surjection, where $X$ is compact and $Y$ is Hausdorff. Then $f$ is a quotient map. The proof follows directly from \cite[Thm.~3.1.12]{En89}.

\begin{remark}
    By the fact quoted above, both maps $p_X$ and $p_Y$ in definition \ref{definition equivalence with cuts and reglues} are quotient maps with finite equivalence classes. Such quotient maps formalize the intuition of gluing together all points in an equivalence class. Taking a right inverse $p_X^{-1}$, we then formalize cutting the space $X /_\sim$ apart. The composition $f := p_Y \circ p_X^{-1}$ then formalizes the procedure of first cutting a space and then re-gluing parts of that space.
    
    We require $p_X$ to be finite-to-one in order to avoid the trivial case where we choose $Z := X \times Y$. Here, any two spaces of the same cardinality would be equivalent. On a geometric note, this finiteness condition prevents us from performing infinitely many cuts through a single point.
\end{remark}

\begin{example}
    We show that the two spaces $X := \mathbb{S}^1$ and $Y := \mathbb{S}^1 \sqcup \mathbb{S}^1$ are equivalent with cuts and subsequent re-glues. We choose $Z$ as the following subspace.
    $$Z := [-1 , 1 ] \times \{ \pm 1 \} \subseteq \mathbb{R}^2$$
    Set $p_X$ to be the quotient map that identifies the two points $( 1 , \pm 1 )$ and the two points $( -1 , \pm 1 )$ in $Z$ and set $p_Y$ to be the quotient map identifying the two points $( \pm 1 , 1 )$ and the two points $( \pm 1 , -1 )$ in $Z$. To conclude, we choose a right inverse $p_X^{-1}$ of $p_X$ whose image contains the two points $( \pm 1 , \pm 1 ) \in Z$.
\end{example}

\begin{proposition}
\label{proposition equivalence with cuts and re-glues is equivalence}
    The equivalence of spaces with cuts and subsequent re-glues is an equivalence relation on the class of all compact regular Hausdorff spaces.
\end{proposition}

\begin{proof}
    $X$ is assumed to be a compact regular Hausdorff space. Therefore, taking $Z := X$ and $p_X \equiv p_Y \equiv p_X^{-1} := id_X$ gives that $X \bumpeq X$ and $\bumpeq$ is reflexive.
    
    Now let $X \bumpeq Y$ where $p_X^{-1}$ is a right inverse of $p_X$. We define \linebreak$p_Y^{-1} := p_X^{-1} \circ f^{-1}$. We have the following equivalent maps.
    \begin{equation*}
        p_Y \circ p_Y^{-1} \equiv p_Y \circ p_X^{-1} \circ f^{-1} \equiv p_Y \circ p_X^{-1} \circ ( p_Y \circ p_X^{-1} )^{-1} \equiv id_Y
    \end{equation*}
    Thus, $p_Y^{-1}$ is a right inverse of $p_Y$, and $\bumpeq$ is symmetric. Now let $X \bumpeq Y$ and $Y \bumpeq W$. By the symmetry of $\bumpeq$, we may assume $Y \bumpeq X$, and we get the following diagram. All maps and spaces are defined as in definition \ref{definition equivalence with cuts and reglues} above.
    \begin{center}
        \begin{tikzcd}
            &
                Z
                \ar[dl, "p_X"']
                \ar[dr, "p_Y"']
                &
                    &
                        Z'
                        \ar[dl, "q_Y"]
                        \ar[dr, "q_W"]
                        \\
            X
            &
                &
                    Y
                    \ar[ll, "f"]
                    \ar[ul, "p_Y^{-1}"'{pos=0.7}, bend right = 30]
                    \ar[ur, "q_Y^{-1}"{pos=0.7}, bend left = 30]
                    \ar[rr, "g"']
                    &
                        &
                            W
        \end{tikzcd}
    \end{center}
    Note that the space $Z \sqcup Z'$ is compact regular Hausdorff. We define the maps $\pi_X : Z \sqcup Z' \to X$ and $\pi_W : Z \sqcup Z' \to W$ as follows.
    \begin{equation*}
        \pi_X : Z \sqcup Z' \to X , \quad 
    x \mapsto 
    \begin{cases}
        p_X ( x ) \quad &\text{if } x \in Z \\
        ( p_X \circ p_Y^{-1} \circ q_Y ) ( x ) \quad &\text{if } x \in Z'
    \end{cases}
    \end{equation*}
    \begin{equation*}
    \pi_W : Z \sqcup Z' \to W , \quad 
    x \mapsto 
    \begin{cases}
        ( q_W \circ q_Y^{-1} \circ p_Y ) ( x ) \quad &\text{if } x \in Z \\
        q_W ( x ) \quad &\text{if } x \in Z'
    \end{cases}
    \end{equation*}
    The two restrictions $\pi_X |_Z$ and $\pi_X |_{Z'}$ are both continuous by definition. The function $\pi_X$ therefore also is continuous. $\pi_X |_Z$ is a surjection, and thus $\pi_X$ also is surjective. The map $\pi_X |_Z$ is finite-to-one by definition, and $\pi_X |_{Z'}$ is the composition of finite-to-one maps and thus itself finite-to-one. Therefore, $\pi_X$ is a finite-to-one continuous surjection. Analogously, the same holds for $\pi_W$. Define $\pi_X^{-1} := p_Y^{-1} \circ f^{-1} : X \to Z \sqcup Z'$. We now have the following diagram. 
    \begin{center}
        \begin{tikzcd}
            &
                Z \sqcup Z'
                \ar[dl, "\pi_X"]
                \ar[dr, "\pi_W"]
                \\
            X
            \ar[ur, "\pi_X^{-1}", bend left = 30]
            \ar[rr, "g \circ f^{-1}"']
            &
                &
                    W
        \end{tikzcd}
    \end{center}
    We first show that $\pi_X^{-1}$ is in fact a right inverse of $\pi_X$. For every $x \in X$ we have the following.
    \begin{equation*}
        ( \pi_X \circ \pi_X^{-1} ) ( x ) = (\pi_X \circ p_Y^{-1} \circ f^{-1} ) ( x ) = (p_X \circ p_Y^{-1} \circ f^{-1} ) ( x ) = ( f \circ f^{-1} ) ( x ) = x
    \end{equation*}
    That is since $im (p_Y^{-1} ) \subseteq Z$ and $\pi_X |_Z \equiv p_X$. Therefore $\pi_X \circ \pi_X^{-1} \equiv id_X$ as required. For every $x \in X$ we further have the following.
    \begin{equation*}
        ( \pi_W \circ \pi_X^{-1} ) ( x ) = (\pi_W \circ p_Y^{-1} \circ f^{-1} ) ( x ) = (q_W \circ p_Y^{-1} \circ f^{-1} ) ( x ) = ( g \circ f^{-1} ) ( x )
    \end{equation*}
    Since $( p_Y^{-1} \circ f^{-1} ) ( x ) \in Z'$ and $\pi_W |_{Z'} \equiv q_W$. We therefore have that \linebreak$\pi_W \circ \pi_X^{-1} \equiv g \circ f^{-1}$ is the composition of two bijections and therefore itself a bijection. This verifies all conditions in \ref{definition equivalence with cuts and reglues} and shows $X \bumpeq W$. Equivalence with cuts and subsequent re-glues is therefore transitive and thus an equivalence relation on the class of all compact regular Hausdorff spaces.
\end{proof}

\begin{proposition}
\label{proposition equivalence with cuts and re-glues implies split homeomorphic}
    If two compact regular Hausdorff spaces $X$ and $Y$ are equivalent with cuts and subsequent re-glues, then $X$ and $Y$ are split homeomorphic.
\end{proposition}

\begin{proof}
    Let $X \bumpeq Y$ be compact regular Hausdorff spaces that are equivalent with cuts and subsequent re-glues. By definition \ref{definition equivalence with cuts and reglues}, there exists a compact space $Z$ and the following maps.
    $$p_X : Z \to X \qquad 
    p_Y : Z \to Y \qquad 
    p_X^{-1} : X \to Z \qquad
    f := p_Y \circ p_X^{-1} : X \to Y$$
    By definition, we have $p_X \circ p_X^{-1} \equiv id_X$. Now define $p_Y^{-1}$ to be the following right inverse of $p_Y$. 
    $$p_Y^{-1} := p_X^{-1} \circ f^{-1} : Y \to Z$$
    We now have the following maps and spaces as summarized in the following diagram.
    \begin{center}
        \begin{tikzcd}
            &
                Z
                \ar[dl, "p_X"]
                \ar[dr, "p_Y"']
                \\
            X
            \ar[ur, "p_X^{-1}", bend left = 30]
            \ar[rr, "{\scriptscriptstyle 1:1}", "f"']
            &
                &
                    Y
                    \ar[ul, "p_Y^{-1}"', bend right = 30]
        \end{tikzcd}
    \end{center}
    Both maps $p_X$ and $p_Y$ are continuous with compact domain and Hausdorff codomain, and therefore they are quotient maps. We may apply proposition \ref{proposition inverse image multifunction of projection from compact to Hausdorff is pre-multisplit} to see that both maps $p_X^{-1}$ and $p_Y^{-1}$ are multi-split continuous. Both $f$ and $f^{-1}$, by definition, then are the composition of a continuous map with a multi-split continuous map and therefore both multi-split continuous. $f$ thus is a split homeomorphism, and we may conclude.
\end{proof}

\begin{proposition}
\label{proposition split homeomorphic implies equivalence with cuts and subsequent re-glues}
    Let $X$ and $Y$ be compact regular Hausdorff spaces that are split homeomorphic. Then $X$ and $Y$ are equivalent with cuts and subsequent re-glues.
\end{proposition}

\begin{proof}
    Let $f : X \to Y$ be a split homeomorphism. We define $Z$ to be the following space.
    $$Z := gr ( f^* ) \subseteq X \times Y$$
    The space $Z$ is closed by corollary \ref{corollary f^* has closed graph} and is a subset of the compact space $X \times Y$. Therefore $Z$ is compact itself. We define the maps $p_X$ and $p_Y$ to be restrictions of the canonical coordinate projections. 
    $$p_X : Z \to X , \quad ( x , y ) \mapsto x 
    \qquad \text{and} \qquad
    p_Y : Z \to Y , \quad ( x , y ) \mapsto y$$
    By applying proposition \ref{Proposition equivalent definition using closure of the graph} to $f$ and to $f^{-1}$, we see that both $p_X$ and $p_Y$ are finite-to-one. Lastly, define the map $p_X^{-1} : X \to Z , \, x \mapsto \big( x , f(x) \big)$. Then $p_X^{-1}$ clearly is a right inverse of $p_X$. We further have $p_Y \circ p_X^{-1} \equiv f$, which is a bijection by assumption. This shows that $X$ and $Y$ are equivalent with cuts and subsequent re-glues.
\end{proof}

\begin{acknowledgements}
    The authors warmly thank asst. prof. Vinay Madhusudanan for numerous insightful discussions. This paper is the result of a 14-week research internship completed by F.M. at the Manipal Institute of Technology (MIT) under the IAESTE A.s.b.l. program (Ref. No. IN-2024-C2701-MU). F.M. gratefully acknowledges the financial support provided by the Manipal Academy of Higher Education, as well as the guidance and assistance from the local IAESTE committee in Manipal and the Department of Mathematics at MIT, Manipal.
\end{acknowledgements}

\begin{contributions}
    F.M. and A.G. jointly developed the main definitions that form the foundation of the work. F.M. worked out most of the statements and proofs independently, while A.G. contributed by suggesting additional statements and filling in specific gaps. F.M. wrote the manuscript, with A.G. providing feedback and corrections. Both authors reviewed and approved the final manuscript.
\end{contributions}

\bibliography{ref}{}
\bibliographystyle{plain}

\end{document}